\documentclass[aos,MSNbibl,seceqn,dvips]{arximspdf}
\usepackage{mathbh}

\doi{10.1214/12-AOS1050} %kopijuoti is PTS
\volume{40}
\issue{6}
\pubyear{2012}
\firstpage{2765}
\lastpage{2797}

\makeatletter
\newcommand{\rrvert}{\vert}
\newcommand{\llvert}{\vert}

\newcommand{\cal}{\mathcal}

\newtheorem{prop}{Proposition}[section]
\newproclaim{assumption}{Assumption}[section]
\newtheorem{coro}[prop]{Corollary}
\newtheorem{theo}[prop]{Theorem}
\newtheorem{lem}[prop]{Lemma}
\newproclaim{rem}{Remark}[section]
\newproclaim{illust}{Illustration}[section]
\newproclaim{illustcont}[illust]{Illustration (continued)}

\renewcommand{\epsilon}{\varepsilon}

\newcommand{\So}{\beta}
\newcommand{\hSo}{\widehat{\So}}
\newcommand{\Sow}{\beta}
\newcommand{\Soclass}{\cF}
\newcommand{\sol}{\beta}
\newcommand{\hsol}{\widehat\sol}
\newcommand{\tsol}{\widetilde\sol}

\newcommand{\Fsw}{\mathcal{F}_b}

\newcommand{\Op}{\Gamma}
\newcommand{\hOp}{\widehat{\Op}}
\newcommand{\Opw}{\gamma}
\newcommand{\Opclass}{\mathcal{G}}

\newcommand{\hgf}{\widehat{g}}

\newcommand{\hw}{\omega}
\newcommand{\Fhw}{\mathcal{F}_\hw}

\newcommand{\bas}{\psi}
\newcommand{\Hspace}{\Hz}

\newcommand{\hOpset}[1]{\Omega_{#1}}
\newcommand{\hOpsetmn}{\hOpset{m,n}}
\newcommand{\Xiset}[1]{\mho_{#1}}
\newcommand{\Xisetmn}{\Xiset{m,n}}
\newcommand{\aset}[1]{\mathcal{A}_{#1}}
\newcommand{\asetn}{\aset{n}}
\newcommand{\bset}[1]{\mathcal{B}_{#1}}
\newcommand{\bsetn}{\bset{n}}
\newcommand{\cset}[1]{\mathcal{C}_{#1}}
\newcommand{\csetn}{\cset{n}}

\newcommand{\eset}[1]{\mathcal{E}_{#1}}
\newcommand{\esetn}{\eset{n}}

\def\pen{\mathrm{pen}}
\def\hpen{\widehat{\mathrm{pen}}}
\def\bias{\mathrm{bias}}
\def\contr{\Psi}

\newcommand{\Me}{\widehat{M}}

\def\Var{\operatorname{\mathbb{V}ar}}
\def\Ex{{\mathbb E}}
\def\argmin{\mathop{\arg\min}}
\def\Diag{\nabla}

\newcommand{\Hz}{{\mathbb H}}
\newcommand{\Nz}{{\mathbb N}}
\newcommand{\Rz}{{\mathbb R}}
\newcommand{\Sz}{{\mathbb S}}

\newcommand{\one}{{\mathbh1}}
\newcommand{\Id}{[\mathrm{Id}]}

\newcommand{\cA}{{\cal A}}
\newcommand{\cB}{{\cal B}}
\newcommand{\cC}{{\cal C}}
\newcommand{\cE}{{\cal E}}
\newcommand{\cF}{{\cal F}}
\newcommand{\cG}{{\cal G}}
\newcommand{\cR}{{\cal R}}
\newcommand{\cT}{{\cal T}}
\newcommand{\cX}{{\cal X}}

\newcommand{\uk}{{\underline{k}}}
\newcommand{\um}{{\underline{m}}}
\newcommand{\whm}{\widehat{m}}
\newcommand{\hsigma}{\widehat{\sigma}}
\newcommand{\htau}{\widehat{\tau}}

\setattribute{abstract}{width}{290pt}

\makeatother

\begin{document}
\begin{frontmatter}

\title{Adaptive functional linear regression\thanksref{T1}}
\runtitle{Adaptive functional linear regression}

\thankstext{T1}{Supported by the IAP research network no. P6/03 of the
Belgian Government (Belgian Science Policy).}

\begin{aug}
\author[A]{\fnms{Fabienne} \snm{Comte}\ead[label=e1]{fabienne.comte@parisdescartes.fr}}
\and
\author[B]{\fnms{Jan} \snm{Johannes}\corref{}\ead[label=e2]{jan.johannes@uclouvain.be}}
\runauthor{F. Comte and J. Johannes}
\affiliation{Universit\'{e} Paris Descartes and Sorbonne Paris Cit\'{e}, and~Universit\'{e}~Catholique de Louvain}
\address[A]{Laboratoire MAP5\\
UMR CNRS 8145\\
45, rue des Saints-P\`{e}res\\
F-75270 Paris Cedex 06\\
France\\
\printead*{e1}}
\address[B]{Institut de Statistique, Biostatistique\\
\quad et Sciences Actuarielles (ISBA)\\
Voie du Roman Pays 20\\
B-1348 Louvain-la-Neuve\\
Belgium\\
\printead{e2}} %adresu isvedimo komanda gale!
\end{aug}

% HISTORY:
\received{\smonth{12} \syear{2011}}
\revised{\smonth{6} \syear{2012}}

% ABSTRACT
%
\begin{abstract}
We consider the estimation of the slope function in functional linear
regression, where scalar responses are modeled in dependence of random
functions. Cardot and Johannes [\textit{J. Multivariate Anal.}
\textbf{101} (2010) 395--408] have shown that a thresholded projection
estimator can attain up to a constant minimax-rates of convergence in a
general framework which allows us to cover the prediction problem with
respect to the mean squared prediction error as well as the estimation
of the slope function and its derivatives. This estimation procedure,
however, requires an optimal choice of a tuning parameter with regard
to certain characteristics of the slope function and the covariance
operator associated with the functional regressor. As this information
is usually inaccessible in practice, we investigate a fully data-driven
choice of the tuning parameter which combines model selection and
Lepski's method. It is inspired by the recent work of Goldenshluger and
Lepski [\textit{Ann. Statist.} \textbf{39} (2011) 1608--1632]. The
tuning parameter is selected as minimizer of a stochastic penalized
contrast function imitating Lepski's method among a random collection
of admissible values. This choice of the tuning parameter depends only
on the data and we show that within the general framework the resulting
data-driven thresholded projection estimator can attain minimax-rates
up to a constant over a variety of classes of slope functions and
covariance operators. The results are illustrated considering different
configurations which cover in particular the prediction problem as well
as the estimation of the slope and its derivatives. A simulation study
shows the reasonable performance of the fully data-driven estimation
procedure.
\end{abstract}

% KEYWORDS

\end{frontmatter}

%s1 #&#
\section{Introduction}\label{secintro}
In functional linear regression the dependence of a real-valued
response $Y$ on the variation of a random function $X$ is studied.
Typically the functional regressor $X$ is assumed
to be square-integrable or more generally to take its values in a
separable Hilbert space $\Hz$ with inner product $\langle\cdot,\cdot
\rangle_\Hz$ and norm
\mbox{$\Vert\cdot\Vert_\Hz$}. Furthermore, we suppose that $Y$ and $X$ are
centered, which simplifies the notations, and that the dependence
between $Y$ and $X$ is linear in the sense that
%
%e1.1 #&#
\begin{equation}
\label{introe1}Y=\langle\sol,X\rangle_\Hz+\sigma \epsilon,\qquad
\sigma>0,
\end{equation}
for some slope function $\sol\in\Hz$ and error term $\epsilon$ with
mean zero and variance one. Assuming an independent and identically
distributed (i.i.d.) sample of $(Y,X)$, the objective of this paper is
the construction of a fully data driven estimation procedure of
the slope function $\sol$ which still can attain minimax-optimal rates
of convergence.

Functional linear models have become very important in a diverse range
of disciplines, including medicine, linguistics, chemometrics as well
as econometrics; see, for instance,~\cite{FerratyVieu2006} and
\cite{RamsaySilverman2005}, for several case studies, or more specific,
\cite{ForniReichlin1998} and~\cite{PredaSaporta2005} for applications
in economics. The main class of estimation procedures of the slope
function studied in the statistical literature is based on principal
components regression; see, for example,
\cite{Bosq2000,CardotFerratySarda1999,CardotMasSarda2007,FrankFriedman1993}
or~\cite{MullerStadtmuller2005} in the context of generalized linear
models. The second important class of estimators relies on minimizing a
penalized least squares criterion which can be seen as generalization
of the ridge regression; cf.~\cite{CardotFerratySarda2003} and
\cite{MarxEilers1999}. More recently an estimator based on dimension
reduction and threshold techniques has been proposed by Cardot and
Johannes~\cite{CardotJohannes2010} which borrows ideas from the inverse
problems community (\cite{EfromovichKoltchinskii2001} and \cite
{HoffmannReiss2008}). It is worth noting that all the proposed
estimation procedures rely on the choice of at least one tuning
parameter, which in turn, crucially influences the attainable accuracy
of the constructed estimator.

It has been shown, for example, in~\cite{CardotJohannes2010}, that the
attainable accuracy of an estimator of the slope $\So$ is essentially
determined by a priori conditions
imposed on both the slope function and the covariance operator $\Op$
associated to the random function $X$ (defined below). These conditions
are usually captured by suitably chosen classes
$\Soclass\subset\Hz$ and $\Opclass$ of slope functions and covariance
operators, respectively. Typically, the class $\cF$ characterizes the
level of smoothness of the slope function, while the class
$\cG$ specifies the decay of the sequence of eigenvalues of $\Op$. For
example,~\cite{CaiZhou2008,CrambesKneipSarda2007} or \cite
{HallHorowitz2007} consider differentiable slope functions and a
polynomial decay of the eigenvalues of $\Op$.
Furthermore, given a weighted norm $\Vert\cdot\Vert_\hw$ and the completion
$\cF_\hw$ of $\Hz$ with respect to
$\Vert\cdot\Vert_\hw$ we shall measure the performance of an
estimator $\hsol$ of
$\sol$ by its maximal $\cF_{\hw}$-risk over a class $\cF\subset\cF_{\hw
}$ of slope functions and a class
$\cG$ of covariance operators, that is,
\[
R_\hw[\hsol;\cF,\cG]:=\sup_{\sol\in\cF}\sup_{\Op\in\cG}\Ex \|
\hsol-\sol\|_{\hw}^2.
\]
This general framework with appropriate choice of the weighted norm
\mbox{$\Vert\cdot\Vert_\hw$} allows us to cover the prediction problem
with respect to
the mean squared prediction error (see, e.g., \cite
{CardotFerratySarda2003} or
\cite{CrambesKneipSarda2007}) and the estimation not only of the slope
function (see, e.g.,~\cite{HallHorowitz2007}) but also of its
derivatives. For a detailed discussion, we refer to \cite
{CardotJohannes2010}. Having these applications in mind the additional
condition $\cF\subset\cF_{\hw}$ only means
that the estimation of a derivative of the slope function necessitates
its existence. Assuming an i.i.d. sample of $(Y,X)$ of size $n$ obeying
model (\ref{introe1}) Cardot and Johannes
\cite{CardotJohannes2010} have derived a lower bound of the maximal
weighted risk, that is,
\[
R_\hw^*[n;\cF,\cG]\leq C \inf_{\hsol}R_\hw[\hsol;
\cF,\cG]
\]
for some finite positive constant $C$ where the infimum is taken over
all possible estimators $\hSo$. Moreover, they have shown that a
thresholded projection estimator
$\hSo_{{m^*_n}}$ in dependence of an optimally
chosen tuning parameter ${m^*_n}\in\Nz$ can attain this lower bound up
to a
constant $C>0$,
\[
R_\hw[\hSo_{{m^*_n}};\cF,\cG]\leq C R_\hw^*[n;\cF,
\cG]
\]
for a variety of classes $\cF$ and $\cG$. In other words, $R_\hw^*[n;\cF
,\cG]$ is the minimax rate of convergence and
$\hSo_{{m^*_n}}$ is minimax-optimal. The optimal choice ${m^*_n}$ of
the tuning parameter, however, follows from a
classical squared-bias-variance compromise and requires an a priori
knowledge about the classes $\cF$ and $\cG$, which is usually
inaccessible in practice.

In this paper we propose a fully data driven method to select a tuning
parameter $\whm$ in such a way that the resulting data-driven estimator
$\hsol_{\whm}$ can still attain the
minimax-rate $R_\hw^*[n;\cF,\cG]$ up to a constant over a variety of
classes $\cF$ and $\cG$. It is interesting to note that, considering a
linear regression model with infinitely
many regressors, Goldenshluger and Tsybakov \cite
{GoldenshlugerTsybakov2001,GoldenshlugerTsybakov2003} propose an
optimal data-driven prediction procedure allowing sharp oracle
inequalities. However, a
straightforward application of their results is not obvious to us
since they assume a priori standardized regressors, which in turn, in
functional linear regression necessitates
the covariance operator $\Op$ to be fully known in advance. In
contrast, given a jointly normally distributed regressor and error
term, Verzelen~\cite{Verzelen2010} establishes sharp
oracle inequalities for the prediction problem in case the covariance
operator is not known in advance. Although, it is worth noting that
considering the mean prediction error as
risk eliminates the ill-posedness of the underlying problem, which in
turn leads to faster minimax rates of convergences of the prediction
error than, for example, the mean
integrated squared error. Cai and Zhou~\cite{CaiZhou2008} present a
fully data-driven estimation procedure of the slope function which
attains optimal rates of convergence with respect to the
maximal mean integrated squared error. On the other hand, covering
both of these two risks within the general framework discussed above,
Comte and Johannes~\cite{ComteJohannes2010} consider functional
linear regression with circular functional regressor which results in
a partial knowledge of the associated covariance operator, that is, its
eigenfunctions are known in advance, but the
eigenvalues have to be estimated. In this situation, Comte and
Johannes~\cite{ComteJohannes2010} have applied successfully a model
selection approach which is inspired by the work of
\cite{BarronBirgeMassart1999} now extensively discussed in \cite
{Massart07}. In the circular case, it is possible to develop the
unknown slope function in the eigenbasis of the
covariance operator, which in turn, allows one to derive an orthogonal
series estimator in dependence of a dimension parameter. This dimension
parameter has been chosen
fully data driven by a model selection approach, and it is shown that
the resulting data-driven orthogonal series
estimator can attain minimax-optimal rates of convergence up to a
constant. Although, the proof crucially relies on the possibility to
write the orthogonal series estimator as a
minimizer of a contrast.

In this paper we do not impose an a priori knowledge of the eigenbasis,
and hence the orthogonal series estimator is no more accessible to us.
Instead, we consider the thresholded projection estimator $\hsol_{m}$
as presented in~\cite{CardotJohannes2010} which we did not succeed to
write as a minimizer of a contrast. Therefore, our selection method
combines model selection and Lepski's method (cf.~\cite{Lepski1990} and
its recent review in~\cite{Mathe2006}) which is inspired by a bandwidth
selection method in kernel density estimation proposed recently in
\cite {GoldenshlugerLepski2011}. Selecting the dimension parameter
$\whm$ as minimizer of a stochastic penalized contrast function
imitating Lepski's method among a random collection of admissible
values, we show that the fully data-driven estimator $\hSo_{\whm}$ can
attain the minimax-rate up to a constant $C>0$, that is,
%
%e1.2 #&#
\begin{equation}
\label{introe2} R_\hw[\hSo_{\whm};\cF,\cG]\leq C\cdot
R_\hw^\star[n;\cF,\cG]
\end{equation}
for a variety of classes $\cF$ and $\cG$. We shall emphasize that in
contrast to the result obtained in~\cite{CaiZhou2008}, we show that the
proposed estimator can
attain minimax-optimal rates without specifying in advance neither that
the slope function belongs to a class of differentiable or
analytic functions nor that the decay of the eigenvalues is polynomial
or exponential. The only price for this flexibility is in term of the
constant $C$ which is
asymptotically not equal to one; that is, the oracle inequality (\ref
{introe2}) is not sharp.

The paper is organized as follows: in Section~\ref{secmet} we briefly
introduce the thresholded projection estimator $\hSo_m$ as proposed in
\cite{CardotJohannes2010}. We present the
data driven method to select the tuning parameter\vspace*{1pt} and prove a first
upper risk-bound for the fully data-driven estimator $\hSo_{\whm}$
which emphasizes the key arguments.
In Section~\ref{secminimax} we review the available minimax theory as
presented in~\cite{CardotJohannes2010}. Within this general framework
we derive upper risk-bounds for the
fully-data driven estimator imposing additional assumptions on the
distribution of the functional regressor $X$ and the error term
$\epsilon$. Namely, we suppose first that $X$ and $\epsilon$ are
Gaussian random variables and second that they satisfy certain moment
conditions. In both cases the proof of the upper risk-bound employs the
key arguments given in Section~\ref{secmet}, while more technical
aspects are
deferred to the \hyperref[app]{Appendix}. The results in this paper are illustrated
considering different configurations of classes $\cF$ and $\cG$. We
recall the minimax-rates in this
situations and show that up to a constant, these rates are attained by
the fully-data driven estimator. A simulation study illustrating the
reasonable performance of the fully
data-driven estimation procedure is available at the supplementary
material archive.

%%%%%%%%%%%
%s2 #&#
\section{Methodology}\label{secmet}
Consider the functional linear model (\ref{introe1}) where the random
function $X$ and the error term $\epsilon$ are independent. Let the
centred random function $X$, that is,
$\Ex\langle X,h\rangle_\Hz=0$ for all $h\in\Hz$, have a finite
second moment,
that is, $\Ex\Vert X\Vert_\Hz^2<\infty$. Multiplying both sides in
(\ref{introe1}) by $\langle X,h\rangle_\Hz$ and taking the
expectation leads to the normal equation
%
%e2.1 #&#
\begin{equation}
\label{mete1}\qquad\langle g,h\rangle_\Hz:=\Ex\bigl[Y\langle X,h
\rangle_\Hz\bigr]= \Ex\bigl[\langle\sol,X\rangle_\Hz
\langle X,h\rangle_\Hz\bigr]=:\langle\Op \sol,h\rangle_\Hz\qquad
\forall h\in\Hz,
\end{equation}
where $g$ belongs to $\Hz$, and $\Op$ denotes the covariance operator
associated to the random function $X$. Throughout the paper we shall
assume that there exists a solution
$\sol\in\Hz$ of equation (\ref{mete1}) and that the covariance
operator $\Op$ is strictly positive definite which ensures the
identifiability of the slope function $\sol$; cf. \cite
{CardotFerratySarda2003}. However, due to the finite second moment of
$X$ the associated covariance operator $\Op$ has a finite trace; that
is, it is nuclear. Thereby,
solving equation (\ref{mete1})
is an \textit{ill-posed inverse problem} with the additional difficulty
that $\Op$ is unknown and has to be
estimated; for a detailed discussion of ill-posed inverse problems in
general we refer to~\cite{EHN00}.

%s2.1 #&#
\subsection{Thresholded projection estimator}
In this paper, we follow~\cite{CardotJohannes2010} and consider a
linear Galerkin approach to derive an estimator of the slope function
$\sol$. Here and subsequently, let $\{\bas_j\}_{j\geq1}$ be a
pre-specified orthonormal basis in $\Hz$ which in general does not
correspond to the eigenbasis of the operator $\Op$ defined in
(\ref{mete1}). With respect to this basis, we consider for all
$h\in\Hspace$ the development $h=\sum_{j=1}^\infty[h]_j\bas_j$
where the sequence\vspace*{-1pt} $([h]_j)_{j\geq1}$ with generic elements
$[h]_j:=\langle h,\bas_j\rangle_\Hz$ is square-summable, that is,
$\Vert h\Vert_\Hz^2=\sum_{j\geq1}[h]_j^2<\infty$. Moreover, given any
strictly positive sequence of weights $(\hw_j)_{j\geq1}$ define the
weighted norm $\Vert h\Vert_\hw^2:= \sum_{j=1}^\infty\hw_j
[h]_j^2$. We
will refer to any sequence as a whole by omitting its index as, for
example, in ``the sequence of weights $\hw$.'' Furthermore, for
$m\geq1$ let $[h]_{\um}:=([h]_1,\ldots,[h]_m)^t$ (where $x^t$ is the
transpose of $x$), and let $\Hz_m$ be the subspace of $\Hz$ spanned by
$\{\bas_1,\ldots,\bas_m\}$. Obviously,\vspace*{1pt} the norm of $h\in\Hz_m$ equals
the Euclidean norm of
its coefficient vector $[h]_{\um}$, that is, $\Vert h\Vert_\Hz
=([h]_\um^t[h]_\um)^{1/2}=:\Vert[h]_\um\Vert$ with a slight abuse of
notation.

An element $\So^{m}\in\Hspace_m$ satisfying
$\Vert g-\Op\So^{m}\Vert_\Hz\leq\Vert g-\Op\breve\sol\Vert_\Hz$
for all
$\breve\sol\in\Hz_m$,
is called a Galerkin solution of equation (\ref{mete1}).
Since the covariance operator $\Op$ is strictly positive definite, it
follows that the covariance matrix $[\Op]_{\um}:=\Ex([X]_{\um
}[X]_{\um}^t)$
associated\vspace*{1pt} with the $m$-dimensional random vector $[X]_{\um}$ is
strictly positive definite too. Consequently, the Galerkin solution\vadjust{\goodbreak}
$\So^{m}\in\Hz_m$ is uniquely determined by
$[\So^{m}]_{\um}=[\Op]_{\um}^{-1}[g]_{\um}$ and $[\So^{m}
]_j=0$ for
all $j>m$. Although, it does generally not correspond to the orthogonal
projection of $\So$ onto the subspace $\Hz_m$ and the approximation
error $\sup_{k\geq m} \Vert\So^{k}-\sol\Vert_\hw$
does generally not converge to zero as $m\to\infty$. Here and
subsequently, however, we restrict ourselves
to classes $\cF$ and $\cG$ of slope functions and covariance operators,
respectively, which ensure the convergence. Obviously, this is a minimal
regularity condition for us since we
aim to estimate the Galerkin solution.

Assuming a sample $\{(Y_i,X_i)\}_{i=1}^n$ of $(Y,X)$ of size $n$, it is
natural to consider the estimators\vspace*{2pt}
$\hgf:=n^{-1}\sum_{i=1}^n Y_i X_i$ and $\hOp:=n^{-1}\sum_{i=1}^n
\langle\cdot, X_i\rangle_\Hz X_i$ for $g$ and $\Op$, respectively. Moreover,
let $[\hOp]_{\um}:=n^{-1}\sum_{i=1}^n [X_i]_{\um}[X_i]_{\um}^t$ and
note that $[\hgf]_{\um}= n^{-1}\sum_{i=1}^n Y_i[X_i]_{\um}$.
Replacing\vspace*{1pt}
the unknown quantities by their empirical counterparts
$\tsol{}^m\in\Hz_m$ denotes a Galerkin solution satisfying
$\Vert\hgf-\hOp\tsol{}^m\Vert_\Hz\leq\Vert\hgf-\hOp\breve
\sol\Vert_\Hz$ for
all $\breve\sol\in\Hz_m$. Observe that there exists always a solution
$\tsol{}^m$, but it might not be unique. Obviously, if $[\hOp]_{\um}$ is
nonsingular, then $[\tsol{}^m]_{\um}=[\hOp]_{\um}^{-1} [\hgf]_{\um
}$. We
shall emphasize the multiplication with the inverse of the random
matrix $[\hOp]_{\um}$ which may result in an unstable estimator even in
case $[\Op]_{\um}$ is well conditioned. Let
$\one_{\{\Vert[\hOp]^{-1}_{\um}\Vert_s\leq n\}}$ denote the indicator
function which takes the value one if $[\hOp]_{\um}$ is nonsingular
with spectral norm\vspace*{2pt}
$\Vert[\hOp]_{\um}^{-1}\Vert_s:=\sup_{\Vert z\Vert=1}\Vert[\hOp
]_{\um}^{-1}z\Vert$
of its inverse bounded by $n$, and the value zero otherwise. The
estimator of $\sol$ proposed in~\cite{CardotJohannes2010} consists of
thresholding the estimated Galerkin solution, that is,
%
%e2.2 #&#
\begin{equation}
\label{introestimator} \hsol_m:= \tsol{}^m
\one_{\{\Vert[\hOp]^{-1}_{\um}\Vert_s\leq n\}}.
\end{equation}
In the next paragraph we introduce a data-driven method to select the
dimension parameter $m\in\Nz$.

%s2.2 #&#
\subsection{Data-driven thresholded projection estimator}
Given a random integer $\Me$ and a random sub sequence of penalties
$(\hpen_1,\ldots,\hpen_{\Me})$, we select the dimension parameter
$\whm
$ among the random collection of admissible values $\{1,\ldots,\Me\}$ as
minimizer of a penalized contrast criterion. To be precise, setting
$\argmin_{m\in A}\{a_m\}:=\min\{m\dvtx  a_m\leq a_{m'}, \forall
m'\in A\}$ for a sequence $(a_m)_{m\geq1}$ with minimal value in
$A\subset\mathbb N$,
we define
%
%e2.3 #&#
\begin{equation}
\label{defwhm} \whm:=\argmin_{1\leq m\leq\Me}{ \lbrace\contr_m+
\hpen_m \rbrace}.
\end{equation}
The data-driven estimator of $\So$ is now given by $\hSo_{\whm}$, and
below we derive an upper bound for its maximal $\cF_{\hw}$-risk. The
choice of the
$\cF_\hw$-risk as performance measure is reflected in the definition of
the contrasts, that is,
\[
\contr_m:=\max_{m\leq k\leq\Me}{ \bigl\lbrace\Vert
\hsol_{k}-\hsol_{m}\Vert^2_\hw-
\hpen_k \bigr\rbrace},\qquad 1\leq m\leq\Me.\vadjust{\goodbreak}
\]
The construction of the random
penalty sequence $\hpen$ and the upper bound $\Me$ given below is
guided by the key arguments used in the proof of the $\cF_\hw$-risk
bound which we present first.
A central step for our reasoning is the next assertion which employs
essentially the
particular choice of the contrast.
%
%le2.1 #&#
\begin{lem}\label{elementarybound}Consider the approximation errors
$\bias_m=\sup_{m\leq k} \Vert\So^k-\So\Vert_\hw$, $m\geq1$.
If the sub sequence $(\hpen_1,\ldots,\hpen_{\Me})$ is nondecreasing,
then we have
%
%e2.4 #&#
\begin{equation}
\label{elementarybounde1}\qquad
\Vert\hsol_{\whm}-\sol\Vert^2_\hw
\leq7\hpen_m +78\bias_m^2+42
\max_{m\leq k\leq\Me} { \biggl( \bigl\|\hsol_{k}-\sol^{k}
\bigr\|^2_\hw -\frac{1}{6}\hpen_k \biggr)
}_+
\end{equation}
for all $1\leq m\leq\Me$, where $(a)_+=\max(a,0)$.
\end{lem}
\begin{pf}
From the definition of $\whm$ we deduce for all $1\leq m\leq\Me$ that
%
%e2.5 #&#
\begin{eqnarray}
\label{prelementarybounde1} \Vert\hsol_{\whm}-\sol
\Vert^2_\hw &\leq&3{ \bigl\lbrace\Psi_m +
\hpen_{\whm} +\Psi_{\whm}+\hpen_m + \Vert
\hsol_{m}-\sol\Vert_\hw^2 \bigr\rbrace}
\nonumber\\[-8pt]\\[-8pt]
&\leq&6{ \lbrace\Psi_m +\hpen_m \rbrace}+3 \|
\hsol_{m}-\sol\|^2_\hw.
\nonumber
\end{eqnarray}
First, employing an elementary triangular inequality allows us to write
\[
\|\hsol_{m}-\sol\|^2_\hw\leq
\frac{1}{3}\hpen_m + 2\bias_m^2+ 2
\max_{m\leq k\leq M} { \biggl( \bigl\|\hsol_{k}-\sol^{k}
\bigr\|^2_\hw -\frac{1}{6}\hpen_k \biggr)
}_+
\]
for all $1\leq m\leq\Me$. Second, since $(\hpen_1,\ldots,\hpen_{\Me})$
is nondecreasing and $4\bias^2_m\geq\max_{m\leq k\leq\Me} \|
\sol^{k}-\sol^{m}\|^2_\hw$, $1\leq m\leq\Me$, it is easily
verified that
\[
\Psi_m\leq6\sup_{m\leq k\leq\Me} { \biggl( \bigl\|\hsol_{k}-
\sol^{k}\bigr\|^2_\hw-\frac{1}{6}
\hpen_k \biggr) }_+ + 12\bias_m^2.
\]
Combining the last two inequalities and (\ref{prelementarybounde1}),
we obtain the result.% assertion \eqref{elementarybounde1}, which
%completes the proof.
\end{pf}
Keeping the last assertion in mind we decompose the $\cF_\hw$-risk with
respect to an event on which the quantities $\hpen_m$ and $\Me$ are
close to some theoretical
counterparts $\pen_m$, $M^-_n$ and $M^+_n$. More precisely, define the event
%
%e2.6 #&#
\begin{equation}
\label{defesetn}\qquad\esetn:= { \bigl\lbrace\pen_k\leq
\hpen_k\leq72\pen_k, \forall 1\leq k\leq M^+_n
\bigr\rbrace }\cap{ \bigl\lbrace M^-_n\leq\Me\leq M^+_n
\bigr\rbrace}
\end{equation}
and the corresponding risk decomposition
%
%e2.7 #&#
\begin{equation}
\label{defdecomp}\qquad R_\hw(\hSo_{\whm};\cF,\cG) =
\sup_{\So\in\cF}\sup_{\Op\in
\cG}\Ex \bigl( \|\hSo_{\whm}-\So
\|^2_\hw\one_{\esetn} \bigr) +
\sup_{\So
\in\cF}\sup_{\Op\in\cG}\Ex \bigl( \|\hSo_{\whm}-\So
\|^2_\hw\one_{\esetn^c} \bigr).\hspace*{-28pt}
\end{equation}
Consider the first right-hand side (r.h.s.) term. If
$(\hpen_1,\ldots,\hpen_{\Me})$ is nondecreasing, then we may apply
Lemma~\ref{elementarybound} which on the event $\esetn$ implies
\begin{eqnarray*}
\|\hSo_{\whm}-\So\|^2_\hw
\one_{\esetn}&\leq&582\max\bigl(\pen_{{m^\diamond_n}}, \bias^2_{{m^\diamond_n}}
\bigr) \\
&&{}+ 42\max_{{m^\diamond_n}\leq k\leq M^+_n} \biggl(\bigl\|\hsol_{k}-
\sol^{k}\bigr\|^2_\hw-\frac{1}{6}
\pen_k \biggr)_+,
\end{eqnarray*}
where ${m^\diamond_n}$ realizes a penalty-squared-bias compromise
among the
collection of admissible values ${ \lbrace1,\ldots,M^-_n
\rbrace}$. Keeping in mind
that ${m^\diamond_n}$ should mimic the value of
the optimal variance-squared-bias trade-off, we wish the upper bound
$M^-_n$ to be as large as possible. In contrast, in order to control the
remainder term, the second
r.h.s. term, we are forced to use a rather small upper bound $M^+_n\geq
M^-_n
$ to ensure that the penalty term is uniformly bounded with increasing
sample size. However, we bound the
remainder term by imposing
the following assumption, which though holds true for a wide range of
classes $\cF$ and $\cG$ under reasonable assumptions on the
distribution of $\epsilon$ and $X$; see Propositions
\ref{gaussp1} and~\ref{momp1} in Section~\ref{secminimax}.
%
%as2.1 #&#
\begin{assumption}\label{assremainder-i} There exists a constant $K_1$
such that
\[
\sup_{\sol\in\cF}\sup_{\Op\in\cG}\Ex{ \biggl\lbrace\max_{{m^\diamond_n}\leq k\leq M^+_n} {
\biggl( \bigl\|\hsol_{k}-\sol^{k}\bigr\|^2_\hw-
\frac{1}{6}\pen_k \biggr) }_+ \biggr\rbrace}\leq
K_1 n^{-1} \qquad\mbox{for all }n\geq1.
\]
\end{assumption}
Roughly speaking, the penalty term $\pen_m$ should provide an upper
bound for the estimator's variation which allows us to establish a
concentration inequality for the
\mbox{$\|\cdot\|_\hw$}-norm of the corresponding empirical process.
However, under
Assumption~\ref{assremainder-i} we bound the first r.h.s. term in
(\ref{defdecomp}) by
%
%e2.8 #&#
\begin{equation}
\label{defdecompe1} \sup_{\So\in\cF}\sup_{\Op\in\cG}\Ex\bigl( \|
\hSo_{\whm}-\So \|^2_\hw\one_{\esetn}
\bigr)\leq582 \sup_{\So\in\cF}\sup_{\Op\in\cG} \max\bigl\{
\pen_{{m^\diamond_n}
},\bias_{{m^\diamond_n}}^2\bigr\}+ 42
\frac{K_1}{n}.\hspace*{-35pt}
\end{equation}
It remains to consider the second r.h.s. term. The conditions on the
distribution of $\epsilon$ and $X$ presented in the next section are
also sufficient to show that the following assumption holds true.
%
%as2.2 #&#
\begin{assumption}\label{assremainder-ii}
There exists a constant $K_2>0$ such that
\[
\sup_{\sol\in\cF}\sup_{\Op\in\cG}\Ex \bigl(\|\hsol_{\whm
}-\sol
\|^2_\hw \one_{\cE^c_n} \bigr)\leq K_2
n^{-1} \qquad\mbox{for all }n\geq1.
\]
\end{assumption}
Under Assumption~\ref{assremainder-ii}, $\Me$ and $\hpen_m$ behave
similarly to their theoretical counterparts with sufficiently high
probability so as not to deteriorate
the estimators risk. The next assertion provides an upper bound for the
maximal $\cF_\hw$-risk over the classes $\cF$ and $\cG$ of the
thresholded projection estimator $\hsol_{\whm}$ with data-driven choice
$\whm$ given by (\ref{defwhm}).
%
%pr2.2 #&#
\begin{prop}\label{methp1} Suppose that $(\hpen_1,\ldots,\hpen_{\Me})$
is nondecreasing. If Assumptions~\ref{assremainder-i} and
\ref{assremainder-ii} hold true, then for all $n\geq1$ we have
$\cR_\hw[\hsol_{\whm};\cF,\cG]\leq582 \sup_{\So\in\cF}\sup_{\Op\in\cG}
\max\{\pen_{{m^\diamond_n}},\bias_{{m^\diamond_n}}^2\}+ (42 K_1+K_2)
n^{-1}$.
\end{prop}
\begin{pf}Keeping in mind the risk decomposition (\ref{defdecomp}) the
upper bound (\ref{defdecompe1}) and Assumption~\ref{assremainder-ii}
imply the result.
\end{pf}
%
%re2.1 #&#
\begin{rem}The first r.h.s. term in the last upper risk-bound is strongly
reminiscent of a variance-squared-bias decomposition of the $\cF_\hw
$-risk for the estimator
$\hsol_{{m^\diamond_n}}$ with dimension parameter ${m^\diamond_n}$.
Indeed, in many
cases the penalty term $\pen_{{m^\diamond_n}}$ is in the same order
as the
variance of the estimator $\hsol_{{m^\diamond_n}}$; cf. Illustration
\ref{illust1}[P-P] and [E-P] below. Consequently,
in this situation the upper risk bound of the data-driven estimator is
essentially given
by $\cR_\hw[\hsol_{{m^\diamond_n}};\cF,\cG]$. Moreover, by
balancing penalty
and squared-bias ${m^\diamond_n}$ just
realizes the optimal trade-off between variance and squared-bias which
in turn in many cases means that $\cR_\hw[\hsol_{{m^\diamond
_n}};\cF,\cG
]$ is
of optimal order.\looseness=1
\end{rem}

We complete this section by introducing our choice for the random upper
bound $\Me$ and the random penalty $\hpen_m$ which takes its
inspiration from~\cite{ComteJohannes2010}. Let us first define some
auxiliary quantities required in the
construction. For $m\geq1$, let $[\Diag_\hw]_{\um}$ denote the
$m$-dimensional diagonal matrix with diagonal entries $(\hw_j)_{1\leq
j\leq m}$, and for any sequence $[K]:=([K]_{\uk})_{k\geq1}$ of
matrices, define\looseness=1
%
%e2.9 #&#
\begin{eqnarray}
\label{defdelta}
\Delta_m^{[K]}&:=&\max_{1\leq k\leq m}\bigl\|[
\Diag_{\hw
}]_{\uk}^{1/2}[K]_{\uk}^{-1}[
\Diag_{\hw}]_{\uk}^{1/2}\bigr\|_{s}
\quad\mbox{and}
\nonumber\\[-8pt]\\[-8pt]
\delta_m^{[K]}&:=&m \Delta_m^{[K]}
\frac{\log(\Delta_m^{[K]}\vee
(m+2))}{\log(m+2)}.
\nonumber
\end{eqnarray}\looseness=0
For $n\geq1$, set
$M^{\hw}_n:=\max{ \lbrace1\leq m\leq{\lfloor n^{1/4}\rfloor}\dvtx
\hw_{(m)}\leq n \rbrace}$ with integer part
${\lfloor n^{1/4}\rfloor}$ of $n^{1/4}$ and $\hw_{(m)}:=\max_{1\leq k\leq
m}\hw_{k}$. For any sequence
$a:=(a_m)_{m\geq1}$ let
%
%e2.10 #&#
\begin{equation}
\label{defM} M_n(a):=\min{ \biggl\lbrace2\leq m\leq
M^{\hw}_n\dvtx  m \hw_{(m)} a_m>
\frac{n}{1+\log n} \biggr\rbrace}-1,
\end{equation}
where we set $M_n(a):=M^{\hw}_n$ if the defining set is empty. Given the
sequence of covariance matrices $[\Op]=([\Op]_\um)_{m\geq1}$
associated with the regressor $X$, define%
%
%e2.11 #&#
\begin{eqnarray}
\label{defpen}\qquad
\pen_m&:=& \kappa \sigma_m^2
\delta_m^{[\Op]} n^{-1} \qquad\mbox{with }
\sigma_m^2:=2 \bigl(\Ex Y^2+[g
]_{\um}^t[\Op]_{\um}^{-1}[g]_{\um}
\bigr) \quad\mbox{and}
\nonumber\\[-8pt]\\[-8pt]
M^\Op&:=& M_n(a) \qquad\mbox{with } a:=\bigl(\bigl\|[
\Op]_{\um}^{-1}\bigr\|_{s}\bigr)_{m\geq1},
\nonumber
\end{eqnarray}
where $\kappa$ is a positive numerical constant to be chosen below.
Roughly\vspace*{1pt} speaking the penalty term provides an upper bound
of the variance of the estimator $\hsol_{m}$ and is in many cases even
in the same order. Its construction, however, allows a deterioration to
ensure that Assumption~\ref{assremainder-i}
can be satisfied; cf. Illustration~\ref{illust1}[P-E].
Moreover, for growing sample size $n$ the penalty sequence is uniformly
bounded over the collection of admissible values
$\{1,\ldots,M_n^\Op\}$. Note that the penalty and the
upper bound still depend on unknown quantities which, however, can
easily be estimated, that is,\looseness=1
%
%e2.12 #&#
\begin{eqnarray}
\label{defhpenMen} \hpen_m&:=&14 \kappa
\hsigma_m^2 \delta_m^{[\widehat{\Op}]}
n^{-1} \nonumber\\
&&\qquad\hspace*{-40pt}\mbox{with } \hsigma_m^2:=2 \Biggl(
\frac
{1}{n}\sum_{i=1}^n
Y^2_i+[\hgf]_{\um}^t[
\hOp]_{\um}^{-1}[\hgf ]_{\um} \Biggr)\quad\mbox{and}\hspace*{-12pt}
\\
\Me&:=& M_n(a) \qquad\mbox{with } a:=\bigl(\bigl\|[\widehat{
\Op}]_{\um
}^{-1}\bigr\|_{s}\bigr)_{m\geq1}.
\nonumber
\end{eqnarray}\looseness=0
Note that by construction $(\hpen_1,\ldots,\hpen_{\Me})$ is
nondecreasing. Indeed, the identity
$\langle\hOp(\hSo_k-\hSo_m),(\hSo_k-\hSo_m)\rangle_\Hz=[\hgf
]_{\uk
}^t[\widehat{\Op}]_{\uk
}^{-1}[\hgf]_{\uk}-[\hgf]_{\um}^t[\widehat{\Op}]_{\um}^{-1}[\hgf
]_{\um}$ holds
true for all $1\leq m\leq
k\leq\Me$. Since $\hOp$ is positive definite, $[\hgf]_{\um
}^t[\widehat{\Op}]_{\um
}^{-1}[\hgf]_{\um}\leq[\hgf]_{\uk}^t[\widehat{\Op}]_{\uk
}^{-1}[\hgf]_{\uk
}$ and
$\hsigma{}^2_m\leq\hsigma_k^2$ which in
turn implies the assertion. Consequently, we may apply Proposition
\ref{methp1} if Assumptions~\ref{assremainder-i} and
\ref{assremainder-ii} hold true.

%%%%%%%%%%%
%s3 #&#
\section{Minimax-optimality}\label{secminimax}
In this section we recall first a general framework proposed by Cardot
and Johannes~\cite{CardotJohannes2010} which allows us to derive
minimax-optimal rates for the maximal $\cF_\hw$-risk,
$\sup_{\sol\in\cF}\sup_{\Op\in\cG}\Ex\|\hsol-\sol\|_{\hw}^2$, over
the classes $\cF$ and~$\cG$.

%s3.1 #&#
\subsection{Notations and basic assumptions}
The classes $\cF$ and $\cG$ of slope functions and covariance
operators, respectively, are
characterized by different weighted norms in $\Hz$ with respect to the
pre-specified orthonormal basis $\{\bas_j,j\in\Nz\}$. Given a strictly
positive sequence of weights $b$ and a radius $r>0$, let $\Fsw$ be
the completion of
$\Hz$ with respect to the weighted norm \mbox{$\|\cdot\|_b$} and the
ellipsoid
$\mathcal{F}^r_b:= \{h\in\Fsw\dvtx
\| h\|_{b}^2\leq r \}$ be the class of possible slope
functions. Furthermore, as usual in the context of ill-posed inverse
problems, we link the mapping properties of the
covariance operator $\Op$ and the regularity condition $\sol\in
\mathcal{F}^r_b$.
Denote by $\mathcal{N}$
the set of all strictly positive nuclear operators defined on $\Hz$.
Given a strictly positive sequence of weights $\Opw$ and a constant
$d\geq 1$ define the class
of covariance operators by%
\[
\cG_{\Opw}^d:= \bigl\{ T\in\mathcal{N}\dvtx  d^{-2}
\| f\|_{\Opw
^2}^2\leq\| T f\|^2
\leq {d^2} \| f\|_{\Opw
^2}^2, \forall f
\in\Hz \bigr\},
\]
where arithmetic operations on sequences are defined element-wise, for
example, $\Opw^{2}=(\Opw^{2}_j)_{j\geq1}$. Let us briefly discuss the
last definition. If $T\in\cG_{\Opw}^d$, then we have $d^{-1}\leq
\langle T\bas_j,\bas_j\rangle/\Opw_j\leq d$, for all $j\geq1$.
Consequently, the
sequence $\Opw$ is necessarily
summable, because $T$ is nuclear. Moreover, if $\lambda$ denotes the
sequence of eigenvalues of $T$, then $d^{-1}\leq\lambda_j/\Opw_j\leq d$, for all $j\geq1$.
In other words the sequence $\Opw$ characterizes the decay of the
eigenvalues of $T\in\cG_{\Opw}^d$. We do not specify the sequences
of weights
$\hw$, $b$ and $\Opw$, but impose from now on the
following minimal regularity conditions.
%
%as3.1 #&#
\begin{assumption}\label{assreg} Let $\hw$,\vspace*{1pt} $b$ and $\Opw$
be strictly positive sequences of weights with $b_1= \hw_1= \Opw_1=
1$, and $\sum_{j=1}^\infty\Opw_j<\infty$ such that the sequences
$b^{-1}$, $\hw b^{-1}$, $\Opw$ and $\Opw^2\hw^{-1}$ are monotonically
nonincreasing
and converging to zero.
\end{assumption}
The last assumption is fairly mild. For example, assuming that $\hw b^{-1}$ is nonincreasing, ensures that $\mathcal{F}^r_b\subset\Fhw$.
Furthermore, it is shown in
\cite{CardotJohannes2010} that the minimax rate $R^*_\hw[n;\mathcal{F}^r_b
,\cG_{\Opw}^d
]$ is of order $n^{-1}$ for all sequences $\Opw$ and $\hw$ such that
$\Opw^2\hw^{-1}$ is nondecreasing.
We will illustrate all our results considering the following three
configurations for the sequences $\hw$, $b$ and $\Opw$.
%
%il3.1 #&#
\begin{illust}\label{illust1}
In all three cases, we take $\hw_j=j^{2s}$, $j\geq1$. Moreover, let:
\begin{longlist}[{[P-E]}]
\item[{[P-P]}] $b_j=j^{2p}$ and $\Opw_j= j^{-2a}$,
$j\geq
1$, with $p> 0$, $a>1/2$ and $p> s>-2a$;
\item[{[E-P]}] $b_j=\exp(j^{2p}-1)$ and $\Opw_j= j^{-2a}$,
$j\geq1$, with $p> 0$, $a>1/2$, $s> -2a$;
\item[{[P-E]}] $b_j=j^{2p}$ and $\Opw_j= \exp(-j^{2a}+1)$,
$j\geq1$, with $p> 0$, $a>0$, and $p> s$;
\end{longlist}
then Assumption~\ref{assreg} is satisfied in all cases.
\end{illust}
%
%re3.1 #&#
\begin{rem} In the configurations [P-P] and
[E-P], the case $s=-a$ can be interpreted as mean-prediction error;
cf.~\cite{CardotJohannes2010}.
Moreover, if ${ \lbrace\bas_j \rbrace}$ is the
trigonometric basis and the value of
$s$ is an integer, then the weighted norm $\| h\|_\hw$ corresponds
to the $L^2$-norm
of the weak $s$th derivative of $h$; cf.~\cite{Neubauer1988}. In other
words in this situation we consider as risk the mean integrated squared
error when estimating the $s$th
derivative of $\sol$. In the configurations [P-P] and
[P-E], the additional condition $p>s$ means that the slope
function has at least $p\geq s+1$
weak derivatives, while for a value $p>1$ in [E-P], the
slope function is assumed to be an analytic function; cf.~\cite{Kawata1972}.
\end{rem}
%
%s3.2 #&#
\subsection{Minimax optimal estimation reviewed}
Let us first recall a lower bound of the maximal $\cF_\hw$-risk over
the classes $\mathcal{F}^r_b$ and $\cG_{\Opw}^d$
due to~\cite{CardotJohannes2010}. Given an i.i.d. sample of $(Y,X)$ of
size $n$ and sequences as in Assumption~\ref{assreg}, define
%
%e3.1 #&#
\begin{eqnarray}
\label{defmstarn} {m^*_n}&:=&\argmin_{m\geq1} \Biggl\{\max
\Biggl(\frac{\hw_m}{b_m}, \sum_{j=1}^m
\frac{\hw_j}{n\Opw_j} \Biggr) \Biggr\} \quad\mbox{and}\nonumber\\[-8pt]\\[-8pt]
{R^*_n}&:=&\max \Biggl(
\frac{\hw_{{m^*_n}}}{b_{{m^*_n}}}, \sum_{j=1}^{{m^*_n}}
\frac{\hw_j}{n\Opw_j} \Biggr).\nonumber
\end{eqnarray}
If $\xi:=\inf_{n\geq1}\{({R^*_n})^{-1}\min(\hw_{{m^*_n}}b_{{m^*_n}
}^{-1},
\sum_{j=1}^{{m^*_n}}{\hw_j}(n\Opw_j)^{-1}) \}>0$, then there exists a
constant $C:=C(\sigma,r,d, \xi)>0$ depending on $\sigma,r,d
$ and $ \xi$ only such that%
%
%e3.2 #&#
\begin{equation}
\label{deflower} \inf_{\tsol} \cR^*_\hw \bigl[\tsol;
\mathcal{F}^r_b,\cG_{\Opw}^d\bigr]
\geq C {R^*_n} \qquad\mbox{for all }n\geq1.
\end{equation}
On the other hand,\vspace*{2pt} considering the dimension parameter ${m^*_n}$ given
in (\ref{defmstarn}) Cardot and Johannes~\cite{CardotJohannes2010}
have shown that the maximal risk
$R^*_\hw[\hsol_{{m^*_n}};\mathcal{F}^r_b,\cG_{\Opw}^d]$ of the
estimator $\hsol_{{m^*_n}
}$ defined in
(\ref{introestimator}) is bounded by ${R^*_n}$ up to constant for a
wide range of sequences $\hw$, $b$ and $\Opw$, provided the random
function $X$ and the error $\epsilon$ satisfy certain additional moment
conditions. In other words ${R^*_n}=R^*_\hw[n;\mathcal{F}^r_b,\cG_{\Opw}^d]$ is the
minimax-rate in this situation, and the estimator $\hsol_{{m^*_n}}$ is
minimax optimal; although, the definition of the dimension parameter
${m^*_n}$ necessitates an a priori knowledge of the sequences $b$ and
$\Opw$. In the remaining part of this paper we show that the
data-driven choice of the dimension parameter constructed in Section
\ref{secmet} can automatically attain the minimax-rate ${R^*_n}$ for a
variety of sequences $\hw$, $b$ and $\Opw$.
First, let us briefly illustrate the minimax result.
%
%il3.2 #&#
\begin{illustcont}\label{illust2} Considering the three configurations
(see Illustration~\ref{illust1}), it has been shown in \cite
{CardotJohannes2010} that the estimator
$\hsol_{{m^*_n}}$ with ${m^*_n}$ as given below attains the rate
$R^*_n$ up to a constant. We write for two strictly positive sequences
$(a_n)_{n\geq1}$ and $(b_n)_{n\geq1}$ that $a_n\sim b_n$, if
$(a_n/b_n)_{n\geq1}$ is bounded away from $0$ and infinity.
\begin{longlist}[{[P-E]}]
\item[{[P-P]}] If $s+a>-1/2$, then $m_n^*\sim
n^{1/(2p+2a+1)}$, while $m_n^*\sim n^{1/[2(p-s)]}$ for $s+a<-1/2$. Thus,
${R^*_n}\sim\max(n^{-(2p-2s)/(2a+2p+1)},n^{-1})$ for $s+a\ne-1/2$. If
$s+a=-1/2$, then $m_n^*\sim(n/\log n)^{1/[2(p-s)]}$ and ${R^*_n}\sim
\log(n)/n$.
\item[{[E-P]}] If $s+a>-1/2$, then ${m^*_n}\sim(\log n
-\frac
{2a+1}{2p}\log(\log n))^{1/(2p)}$ and ${R^*_n}\sim n^{-1} (\log
n)^{(2a+1+2s)/(2p)}$, while
$m_n^*\sim(\log n+(s/p)\log(\log n))^{1/(2p)}$ and ${R^*_n}\sim n^{-1}$
for $s+a<-1/2$ [and ${R^*_n}\sim\log(\log n)/n$ for $a+s=-1/2$].
\item[{[P-E]}] ${m^*_n}\sim(\log n -\frac
{2p+(2a-1)_{+}}{2a}\log(\log n))^{1/(2a)}$ and ${R^*_n}\sim(\log
n)^{-(p-s)/a}$.
\end{longlist}
An increasing value of the parameter $a$ leads in all three cases to a
slower rate ${R^*_n}$, and hence it is called degree of ill-posedness;
cf.~\cite{Natterer84}.
\end{illustcont}
%
%s3.3 #&#
\subsection{Minimax-optimality of the data-driven estimation procedure}
Consider the thresholded projection estimator $\hsol_{\whm}$ with
data-driven choice $\whm$ of the dimension parameter. Supposing that
the joint distribution of the random
function $X$ and the error term $\epsilon$ satisfies certain additional
conditions, we will prove below that Assumptions
\ref{assremainder-i} and~\ref{assremainder-ii}
formulated in Section~\ref{secmet} hold true. These assumptions rely
on the existence of sequences $({m^\diamond_n})_{n\geq1}$ and
$(M^+_n)_{n\geq1}$
which amongst others we define now
referring only to the classes $\mathcal{F}^r_b$ and $\cG_{\Opw}^d$.
Keep\vspace*{-1pt} in mind the
notation given in (\ref{defdelta}) and (\ref{defM}). For $m,n\geq
1$ and
$[\Diag_\Opw]=([\Diag_\Opw]_{\um})_{m\geq1}$ define $\Delta_m^\Opw
:=\Delta_m^{[\Diag_\Opw]}$ and $\delta_m^\Opw:=\delta_m^{[\Diag
_\Opw]}$,
set $M^-_n:= M_{n}(16d^3\Opw^{-1})$ and $M^+_n:= M_{n} ((4d
\Opw
)^{-1} )$, and let
\[
{m^\diamond_n}:=\argmin_{1\leq m\leq M^-_n}{ \biggl\lbrace \max{
\biggl( \frac{\hw_m}{b_m}, \frac{\delta^\Opw_m}{n} \biggr) } \biggr\rbrace}
\quad\mbox{and}\quad
{R^\diamond_n}:=\max{ \biggl( \frac{\hw_{{m^\diamond
_n}}}{b_{{m^\diamond_n}}},
\frac{\delta^\Opw_{{m^\diamond
_n}}}{n} \biggr) },
\]
where ${m^\diamond_n}\leq M^-_n\leq M^+_n$. Let $\Sigma:=\Sigma(\cG_{\Opw}^d)$
denote a
finite constant such that
%
%e3.3 #&#
\begin{equation}
\label{defSigma}\qquad \Sigma\geq\sum_{j\geq1}
\Opw_j \quad\mbox{and}\quad \Sigma\geq \sum_{m\geq1}
\Delta_m^\Opw\exp \biggl(-\frac{m\log(\Delta^\Opw_m \vee
(m+2))}{16(1+\log d)\log(m+2)} \biggr),
\end{equation}
which by construction always exists and depends on the class $\cG_{\Opw}^d$
only. We illustrate below the last definitions by revisiting the three
configurations for the sequences $\hw$, $b$ and $\Opw$ (Illustration
\ref{illust1}).\vspace*{-2pt}
%
%il3.3 #&#
\begin{illustcont}\label{illust3}In the following we state the order
of $M^-_n$ and $\delta_m^\Opw$ which in turn are used to derive the
order of ${m^\diamond_n}$ and ${R^\diamond_n}$.
\begin{longlist}[{[P-E]}]
\item[{[P-P]}] $M^-_n\sim (\frac{n}{1+\log n}
)^{1+2a+(2s)_+} $, $\delta_m^\Opw\sim m^{1+(2a+2s)_+}$ and for
$p>(s)_+$ it follows ${m^\diamond_n}\sim
m^{1/[1+2p-2s+(2a+2s)_+]}$ and ${R^\diamond_n}\sim
n^{-2(p-s)/[1+2p-2s+(2a+2s)_+]}$;
\item[{[E-P]}] $M^-_n\sim (\frac{n}{1+\log n}
)^{1+2a+(2s)_+} $, $\delta_m^\Opw\sim m^{1+(2a+2s)_+}$ and for $p>0$,
${m^\diamond_n}\sim(\log n
-\frac{1+2(a+s)_+-2s}{2p}\log(\log n))^{1/(2p)}$ and ${R^\diamond_n}\sim
n^{-1} (\log n)^{[1+2(a+s)_+]/(2p)}$;
\item[{[P-E]}] $M^-_n\sim(\log n -\frac
{1+2a+2(s)_+}{2a}\log
(\log n))^{1/(2a)}$, $\delta^\Opw_m\sim m^{1+2s+2a}\exp(m^{2a})$ and
for $p>(s)_+$, it follows
${m^\diamond_n}\sim(\log n -\frac{1+2a+2p}{2a}\log(\log
n))^{1/(2a)}$ and
${R^\diamond_n}\sim(\log n)^{-(p-s)/a}$.\vspace*{-2pt}
\end{longlist}
\end{illustcont}
We proceed by formalizing additional conditions on the joint
distribution of $\epsilon$ and $X$, allowing us to prove that
Assumptions~\ref{assremainder-i} and~\ref{assremainder-ii} hold true.\vspace*{-2pt}

\subsubsection*{Imposing a joint normal distribution} Let us first assume
that $X$ is a centred Gaussian $\Hz$-valued random variable; that is,
for all $k\geq1$ and for all finite collections $\{h_1,\ldots,h_k\}
\subset\Hz$ the joint distribution of the real
valued random variables $\langle X,h_1\rangle_\Hz,\ldots, \langle
X,h_k\rangle_\Hz$ is
Gaussian with zero mean vector and covariance matrix with generic elements
$\Ex\langle h_j,X\rangle_\Hz\langle X,h_l\rangle_\Hz$, $1\leq
j,l\leq k$. Moreover,
suppose that the error term is standard normally distributed.
%
%as3.2 #&#
\begin{assumption}\label{assgauss}The joint distribution of $X$ and
$\epsilon$ is normal.\vspace*{-2pt}
\end{assumption}
The more involved proof of the next assertion is deferred to Appendix
\ref{appgauss}.\vspace*{-2pt}%

%pr3.1 #&#
\begin{prop}\label{gaussp1}Assume an i.i.d. $n$-sample of $(Y,X)$
obeying (\ref{introe1}) and Assumption~\ref{assgauss}. Consider
sequences $\hw$, $b$ and $\Opw$ satisfying
Assumption~\ref{assreg} and set $\kappa=96$ in the definition (\ref
{defpen}) and (\ref{defhpenMen}) of the penalty $\pen$ and $\hpen$,
respectively. For the classes $\mathcal{F}^r_b$ and
$\cG_{\Opw}^d$ there exist finite constants $C_1:=C_1(d)$ and
$C_2:=C_2(d)$ depending on $d$ only such that Assumptions
\ref{assremainder-i} and~\ref{assremainder-ii}
with $K_1:=C_1 (\sigma^2+r) \Sigma$ and $K_2:=C_2 (\sigma^2+r)
\Sigma$, respectively, holds true.\vspace*{-2pt}
\end{prop}
By taking the value $\kappa=96$ the random penalty and upper bound
given in (\ref{defhpenMen}) depend indeed only on the data, and hence
the choice $\whm$ in (\ref{defwhm}) is fully data-driven. Moreover, we
can apply Proposition~\ref{methp1}
to
prove the next upper risk-bound for the data-driven thresholded
projection estimator $\hsol_{\whm}$.\vadjust{\goodbreak}

%th3.2 #&#
\begin{theo}\label{gaussth1}Let the assumptions of Proposition \ref
{gaussp1} be satisfied. There exists a finite constant $K:=K(d)$
depending on $d$ only such that
\[
R_\hw\bigl[\hsol_{\whm},\mathcal{F}^r_b,
\cG_{\Opw}^d\bigr]\leq K \bigl(\sigma^2+r\bigr) {
\bigl\lbrace{R^\diamond_n}+ \Sigma n^{-1} \bigr
\rbrace } \qquad\mbox{for all }n\geq1.
\]
\end{theo}
\begin{pf}
We shall provide in the \hyperref[app]{Appendix} among others, the two technical Lemmas
\ref{appprel1} and~\ref{appprel2} which are used in the following.
Moreover, we denote by $K:=K(d)$ a constant depending on $d$ only
which changes from line to
line. Making use of Proposition~\ref{gaussp1}
we intend to apply Proposition~\ref{methp1}.
To this end, if $\sol\in\mathcal{F}^r_b$ and $\Op\in\cG_{\Opw
}^d$, then first
from \ref
{appprel1e4} in Lemma~\ref{appprel1} it follows that $\bias^2_{m^\diamond_n}\leq34
d^8   r  \hw_{m^\diamond_n} b_{m^\diamond_n}^{-1}$ because $\Opw^2\hw^{-1}$
and $\hw b^{-1}$ are nonincreasing due to Assumption~\ref{assreg}.
Second, by combination of
\ref{appprel2e1} and~\ref{appprel2e4} in Lemma~\ref{appprel2},
it is easily verified that $\pen_{m^\diamond}\leq K  (\sigma^2+r)
\delta_{m^\diamond}^\Opw n^{-1} $. Consequently, $\sup_{\sol\in\mathcal
{F}^r_b}\sup_{\Op\in
\cG_{\Opw}^d}\max(\pen_{{m^\diamond_n}},\bias_{{m^\diamond
_n}}^2)\leq K (\sigma^2+r
) {R^\diamond_n}
$ for all $n\geq1$
by combination of the last two estimates and the definition of
${R^\diamond_n}$
which in turn together with Proposition~\ref{methp1} implies the
assertion of the theorem.
\end{pf}

\subsubsection*{Imposing moment conditions} We now dismiss Assumption
\ref{assgauss} and formalize in its place,
conditions on the moments of the random function $X$ and the error
term~$\epsilon$. In particular we use that for all $h\in\Hz$ with
$\langle\Op h,h\rangle=1$, the random variable $\langle h,X\rangle$ is
standardized, that is, has mean zero and variance one.
%
%as3.3 #&#
\begin{assumption}\label{assmom}There exist a finite integer $k\geq
16$ and a finite constant \mbox{$\eta\geq1$} such that $\Ex|\epsilon
|^{4k}\leq\eta^{4k}$ and that for all $h\in\Hz$ with $\langle\Op
h,h\rangle=1$
the standardized random variable $\langle h,X\rangle$ satisfies $\Ex
|\langle h,X\rangle|^{4k}\leq\eta^{4k}$.
\end{assumption}
It is worth noting that for any Gaussian random function $X$ with
finite second moment, Assumption~\ref{assmom} holds true, since for all
$h\in\Hz$ with $\langle\Op h,h\rangle=1$ the random variable $\langle
h,X\rangle$ is standard normally distributed and hence $\Ex|\langle
h,X\rangle|^{2k}=(2k-1)\cdot\cdots\cdot5\cdot 3\cdot 1$. The proof of
the next assertion is again rather involved and deferred to Appendix
\ref{appmom}. It follows, however, along the general lines of the proof
of Proposition~\ref{methp1} though it is not a straightforward
extension. Take as an example the concentration inequality for the
random variable $\|[\Op]^{1/2}_\um([\hgf]_\um- [\hOp]_\um
[\So^{m}]_\um)\|$ in Lemma~\ref{appgaussl2} in Appendix
\ref{appgauss} which due to Assumption~\ref{assgauss} is shown by
employing elementary inequalities for Gaussian random variables. In
contrast, the proof of an analogous result under Assumption \ref
{assmom} given in Lemma~\ref{appgenerall2} in Appendix~\ref{appmom} is
based on an inequality due to Talagrand~\cite{Talagrand1996}
(Proposition \ref {appgeneraltala} in the \hyperref[app]{Appendix} states a version as
presented in~\cite{KleinRio}).
%
%pr3.3 #&#
\begin{prop}\label{momp1}Assume an i.i.d. $n$-sample of $(Y,X)$
obeying (\ref{introe1}) and Assumption~\ref{assmom}. Consider
sequences as in
Assumption~\ref{assreg} and set $\kappa=288$ in the definition
(\ref{defpen}) and (\ref{defhpenMen}) of the penalty $\pen$ and
$\hpen$,
respectively.
For the classes\vadjust{\goodbreak} $\mathcal{F}^r_b$ and $\cG_{\Opw}^d$, there exist
finite constants
$C_1:=C_1(\sigma,\eta,\mathcal{F}^r_b,\cG_{\Opw}^d)$ depending on
$\sigma$, $\eta
$ and
the classes $\mathcal{F}^r_b$ and $\cG_{\Opw}^d$ only, and
$C_2:=C_2(d)$ depending
on $d$ only, such that
Assumptions~\ref{assremainder-i} and~\ref{assremainder-ii} with
$K_1:=C_1 \eta^{64} (\sigma^2+r) \Sigma$ and $K_2:=C_2 \eta^{64}
(\sigma^2+r) \Sigma$, respectively, hold true.
\end{prop}
We remark on a change only in the constants when comparing the last
proposition with Proposition~\ref{gaussp1}. Note further that we need
a larger value for the constant
$\kappa$ than in Proposition~\ref{gaussp1} although it is still a
numerical constant and hence the choice $\whm$ given by (\ref{defwhm})
is again fully
data-driven. Moreover, both values for the constant~$\kappa$, though
convenient for deriving the theory, are far too large in practice. In
our simulation study they are instead determined by means of
preliminary simulations as proposed in
\cite{ComteRozenholcTaupin2006}, for example.
The next\vspace*{1pt} assertion provides an upper risk-bound for the data-driven
thresholded projection estimator $\hsol_{\whm}$ when imposing moment
conditions.
%
%th3.4 #&#
\begin{theo}\label{momth1}Let the assumptions of Proposition \ref
{momp1} be satisfied. There exist finite constants $K:=K(d)$
depending on $d$ only and $K':=K'(\sigma,\eta,\mathcal{F}^r_b,\cG_{\Opw}^d)$
depending on $\sigma$, $\eta$ and the classes $\mathcal{F}^r_b$ and
$\cG_{\Opw}^d$ only
such that
\[
R_\hw\bigl[\hsol_{\whm},\mathcal{F}^r_b,
\cG_{\Opw}^d\bigr]\leq K \bigl(\sigma^2+r\bigr) {
\bigl\lbrace{R^\diamond_n}+ K'
\eta^{64} \Sigma n^{-1} \bigr\rbrace} \qquad\mbox{for all }n\geq1.
\]
\end{theo}
\begin{pf}
Taking into account Proposition~\ref{momp1} rather than Proposition~\ref{gaussp1} we follow line by
line the proof of Theorem~\ref{gaussth1} and we omit the details.
\end{pf}

\subsubsection*{Minimax-optimality}
A comparison of the upper bounds in both Theorems~\ref{gaussth1} and
\ref{momth1} with the lower bound displayed in (\ref {deflower}) shows
that the data-driven estimator $\hsol_{\whm}$ attains\vspace*{-1pt} up
to a constant the minimax-rate ${R^*_n} =\min_{1\leq
m<\infty}\{\max(\frac{\hw_m}{b_m},\sum_{j=1}^m\frac
{\hw_j}{n\Opw_j})\}$ only if ${R^\diamond_n}=\min_{1\leq m\leq
M^-_n}\{\max(\frac{\hw_m}{b_m}, \frac {\delta_m^\Opw}{n})\}$ has the
same order as ${R^*_n}$. Note that, by construction,
$\delta_m^\Opw\geq\sum_{j=1}^m\frac{\hw_j}{\Opw_j}$ for all $m\geq1$.
The next assertion is an immediate consequence of Theorems
\ref{gaussth1} and~\ref{momth1}, and we omit its proof.
%
%co3.5 #&#
\begin{coro}\label{adapcoro} Let the assumptions of either
Theorems \ref
{gaussth1} or~\ref{momth1} be satisfied. If $\xi^\diamond
:=\sup_{n\geq1}\{{R^\diamond_n}
/{R^*_n}\}<\infty$ holds true, then
$R_\hw[\hsol_{\whm};\mathcal{F}^r_b,\cG_{\Opw}^d]\leq C \cdot
\inf_{\tsol}R_\hw
[\tsol;\mathcal{F}^r_b,\cG_{\Opw}^d]$
for all $n\geq1$ and a finite positive constant $C$, where the infimum
is taken over all possible estimators $\tsol$.
\end{coro}
%
%re3.2 #&#
\begin{rem}In the last assertion $\xi^\diamond=\sup_{n\geq1}\{
{R^\diamond_n} /{R^*_n}\}<\infty$ is, for example, satisfied if the
following two conditions hold simultaneously true: (i)~${m^*_n}\leq
M^-_n$ for all $n\geq1$ and (ii) $\Delta_m^\Opw=\max_{1\leq j\leq
m}\hw_j\Opw_j^{-1} \leq Cm^{-1} \sum_{j=1}^m\hw_j\Opw_j^{-1}$ and
$\log(\Delta_m^\Opw\vee (m+2))\leq C\log(m+2)$ for all $m\geq1$.
Observe\vspace*{2pt} that (ii) which implies $\delta_m^\Opw\leq C
\sum_{j=1}^m\frac{\hw_j}{\Opw_j}$ is satisfied in case $\Delta_m^\Opw$
is in the order of a power of $m$ (e.g.,
Illustration~\ref{illust2}[P-P] and [E-P]). If this term has an
exponential order with respect to $m$ (e.g., Illustration
\ref{illust2}[P-E]), then\vspace*{1pt} a deterioration\vadjust{\goodbreak} of the term $\delta_m^\Opw$
compared to the variance term $\sum_{j=1}^m\frac{\hw_j}{\Opw_j}$ is
possible. However,\vspace*{-1pt} no loss in terms of the rate may occur, that is,
$\xi^\diamond<\infty$, when the squared-bias term
$\hw_{m^\diamond_n}b_{m^\diamond_n}^{-1}$ dominates the variance term
$n^{-1}\delta_{m^\diamond_n}^\Opw$; for a detailed discussion in a
deconvolution context, we refer to~\cite{BT1,BT2}.
\end{rem}
Let us illustrate the performance of the data-driven thresholded
projection estimator revisiting the three configurations presented in
Illustration
\ref{illust1}.
%
%pr3.6 #&#
\begin{prop}\label{rate}Assume an i.i.d. $n$-sample of $(Y,X)$
satisfying (\ref{introe1}) and let either Assumptions~\ref{assgauss}
or~\ref{assmom} hold true where we set, respectively,
$\kappa=96$ or $\kappa=288$ in (\ref{defhpenMen}). The data-driven
estimator $\widehat\beta_{{\whm}}$ attains
the minimax-rates ${R^*_n}$, up to a constant, in the three cases
given in
Illustration~\ref{illust1}, if we additionally assume $a+s\geq0$ in
the cases \textup{[P-P]} and \textup{[E-P]}.
\end{prop}
\begin{pf}
Under the stated conditions it is easily verified that the assumptions
of either Theorems~\ref{gaussth1} or~\ref{momth1} are
satisfied. Moreover, the rates
${R^*_n}$ (Illustration~\ref{illust2}) and ${R^\diamond_n}$ (Illustration
\ref
{illust3}) are of the same order if we additionally assume $a+s\geq0$
in the cases [P-P] and [E-P]. Therefore, Corollary~\ref
{adapcoro} applies, and we obtain the assertion.
\end{pf}
%

%sA #&#
\begin{appendix}\label{app}
\begin{center}
APPENDIX\vspace*{10pt}
\end{center}

This section gathers preliminary technical results and the proofs of
Propositions~\ref{gaussp1} and~\ref{momp1}.

%
%sA #&#
\section{Notation}%
We begin by defining and recalling notation to be used in all proofs.
Given $m\geq1$, $\Hz_m$ denotes the subspace of
$\Hz$ spanned by the functions $\{\bas_1,\ldots,\bas_m\}$. $\Pi_m$ and
$\Pi_m^\perp$ denote the orthogonal projections on $\Hz_m$ and its
orthogonal complement $\Hz_m^\perp$,
respectively. If $K$ is an operator mapping $\Hz$ to itself and if we
restrict $\Pi_m K \Pi_m$ to an operator from $\Hz_m$ to itself, then it
can be represented by a matrix $[K]_{\um}$ with generic entries
$\langle\bas_{j},K\bas_{l}\rangle_\Hz=:[K]_{j,l}$ for $1\leq
j,l\leq m$. The
spectral norm of $[K]_{\um} $ is denoted by $\|[K]_{\um}\|_s$, and
the inverse matrix of $[K]_{\um}$
by $[K]_{\um}^{-1}$. Furthermore, $[\Diag_\hw]_{\um}$ and $\Id_{\um}$
denote, respectively, the $m$-dimensional diagonal matrix with diagonal
entries $(\hw_j)_{1\leq j\leq m}$ and the identity
matrix. For $h\in\Hz_m$ it follows $\| h\|_\hw^2= [h]^t_\um
[\Diag_\hw
]_{\um}[h]_\um=\|[\Diag_\hw]_{\um}^{1/2}[h]_\um\|^2$.
Keeping\vspace*{-1pt} in
mind the notation given in
(\ref{defdelta})--(\ref{defhpenMen}) we use\vspace*{-1pt} for all $m\geq1$ in
addition $\Lambda_m^{[\Op]}:=\frac{\log(\Delta_m^{[\Op]}\vee
(m+2))}{\log(m+2)}$,
$\Lambda_m^\Opw:=\frac{\log(\Delta_m^\Opw\vee(m+2))}{\log
(m+2)}$ and
$\Lambda_m^{[\widehat{\Op}]}:=\frac{\log(\Delta_m^{[\widehat{\Op
}]}\vee(m+2))}{\log(m+2)}$
allowing us to write $\delta_m^{[\Op]}= m\Delta_m^{[\Op]}\Lambda_m^{[\Op]
}$, $\delta_m^\Opw=m\Delta_m^\Opw\Lambda_m^\Opw$ and $\delta_m^{[\widehat{\Op}]}=
m\Delta_m^{[\widehat{\Op}]}\Lambda_m^{[\widehat{\Op}]}$.
Given a Galerkin solution $\So^{m}\in\Hz_m$ of equation (\ref{introe2}),
let $Z_m:=Y- \langle\So^{m},X\rangle_\Hz=\sigma\epsilon+\langle
\sol-\So^{m},X\rangle_\Hz$ and denote $\rho_m^2:=\Ex Z^2_m=\sigma^2+\langle\Op(\sol-\So^{m}),
(\sol-\So^{m})\rangle_\Hz$,
$\sigma_Y^2:=\Ex Y^2= \sigma^2+\langle\Op\sol,\sol\rangle_\Hz$ and\vadjust{\goodbreak}
$\sigma_m^2= 2 (\sigma_Y^2+[g]_{\um}^t[\Op]_{\um}^{-1}[g]_{\um} )$
employing that $\epsilon$ and $X$ are uncorrelated.
Define the matrix $[\Xi]_{\um}:= [\Op]_{\um}^{-1/2}[\hOp]_{\um
}[\Op
]_{\um}^{-1/2} - \Id_{\um}$ and the vector $[W]_{\um}:= [\hgf
]_{\um}-
[\hOp]_{\um}
[\So^{m}]_{\um}$
satisfying $\Ex[\Xi]_{\um}= 0$
and $\Ex[W]_{\um} =[\Op(\sol-\So^{m})]_{\um}=0$. Let further
$\hsigma_Y^2:=n^{-1}\sum_{i=1}^nY_i^2$ and define
the events
%
%eA.1 #&#
\begin{eqnarray}
\label{apppree4} \hOpsetmn&:=&\bigl\{ \bigl\|[\hOp]^{-1}_{\um}
\bigr\|_{s}\leq n \bigr\},\qquad \Xisetmn:= \bigl\{ 8\bigl\|[
\Xi]_{\um}\bigr\|_{s}\leq1\bigr\},
\nonumber
\\
\asetn&:=&\bigl\{{1}/{2}\leq\hsigma_Y^2/
\sigma_Y^2\leq{3}/{2}\bigr\},\qquad \bsetn:=\bigl\{\bigl\|[
\Xi]_\uk\bigr\|_{s}\leq1/8,\forall1\leq k\leq
M^{\hw}_n \bigr\},\hspace*{-26pt}
\\
\csetn&:=&\bigl\{8[W]_\uk^t[\Op]_{\uk}^{-1}[W]_\uk
\leq \bigl([g]_\uk^t[\Op]^{-1}_\uk[g]_\uk+
\sigma_Y^2\bigr),\forall 1\leq k\leq
M^{\hw}_n\bigr\},
\nonumber
\end{eqnarray}
and their complements $\hOpsetmn^c$, $\Xisetmn^c$, $\asetn^c$,
$\bsetn^c$ and $\csetn^c$, respectively. Furthermore, we will denote by $C$
universal numerical constants and
by $C(\cdot)$ constants depending only on the arguments. In both cases,
the values of the constants may change from line to line.

%sB #&#
\section{Preliminary results}\label{prpre}
This section gathers results exploiting
Assumption~\ref{assreg} only. The proof of the next lemma can be found
in~\cite{JohannesSchenk2010}.
%
%leB.1 #&#
\begin{lem}\label{appprel1}
Let $\Op\in\cG_{\Opw}^d$ with sequence $\Opw$ as in Assumption
\ref{assreg}.
Then we have:
\renewcommand\thelonglist{(\roman{longlist})}
\renewcommand\labellonglist{\thelonglist}
\begin{longlist}
\item\label{appprel1e1}$\sup_{m\geq1} \{ \Opw_m \|[\Op
]_{\um}^{-1}\|_s \}\leq4d^3$;\vspace*{1pt}
\item\label{appprel1e2}$\sup_{m\geq1}\|[\Diag_{\Opw
}]^{1/2}_{\um}[\Op]_{\um}^{-1}[\Diag_{\Opw}]^{1/2}_{\um}\|_s\leq
4d^3$;\vspace*{1pt}
\item\label{appprel1e3}$\sup_{m\geq1}\|[\Diag_\Opw
]^{-1/2}_{\um}[\Op]_{\um}[\Diag_\Opw]^{-1/2}_{\um}\|_s\leq d$.
\end{longlist}
Let in addition $\sol\in\cF_b^r$ with sequence $b$ as in
Assumption~\ref{assreg}. If $\So^{m}$ denotes a Galerkin solution of
$g=\Op\sol$, then for each strictly positive sequence $w$ such
that $wb^{-1}$ is nonincreasing and for all $m\geq1$ we obtain:

\renewcommand\thelonglist{(\roman{longlist})}
\renewcommand\labellonglist{\thelonglist}
\begin{longlist}
\setcounter{longlist}{3}
\item\label{appprel1e4}$\|\sol-\So^{m}\|_w^2 \leq34
d^8
r  {w_m}{b_m^{-1}}\max ( 1, {\Opw_m^2}{w_m^{-1}} \max_{1\leq
j\leq m} {w_j}{\Opw_j^{-2}} )$;
\item\label{appprel1e5}$\|\So^{m}\|^2_b\leq34  d^8
r$ and $\|\Op^{1/2}(\sol-\So^{m} )\|_\Hz^2\leq34
d^9   r
\Opw_mb_m^{-1}$.
\end{longlist}
\end{lem}
%
%leB.2 #&#
\begin{lem}\label{appprel2} Let Assumption~\ref{assreg} be
satisfied. If $\Op\in\cG_{\Opw}^d$ and $D:=4d^3$, then:%

\renewcommand\thelonglist{(\roman{longlist})}
\renewcommand\labellonglist{\thelonglist}
\begin{longlist}
\item\label{appprel2e1}$d^{-1}\leq\Opw_m\|[\Op]_\um^{-1}\|_s \leq D$,
$d^{-1} \leq\Delta_m^{[\Op]}/\Delta^\Opw_m\leq D $, $(1+\log
d)^{-1}\leq\Lambda_m^{[\Op]}/\break\Lambda^\Opw_m\leq (1+\log D)$,
and $d^{-1}(1+\log d)^{-1}\leq\delta_m^{[\Op]}/\delta^\Opw_m\leq
D(1+\log D)$, for all $m\geq1$;\vspace*{1pt}
\item\label{appprel2e2}$\delta_{M^+_n}^\Opw\leq n 4D(1+\log D)$ and
$\delta_{M^+_n}^{[\Op]} \leq n 4D^2(1+2\log D)$, for all
$n\geq1$;\vspace*{1pt}
\item\label{appprel2e3}$n\geq2\max_{1\leq m\leq M^+_n}\|
[\Op]^{-1}_\um\|$ if $n\geq2D$ and $\hw_{(M^+_n)}M^+_n(1+\log
n)\geq8D^2$;\vspace*{1pt}
\item\label{appprel2e4}$\rho_m^2\leq\sigma_m^2\leq2(\sigma^2+35
d^9r)$, for all $m\geq1$, assuming in addition \mbox{$\sol\in\cF_b^r$}.\vadjust{\goodbreak}
\end{longlist}
\end{lem}
\begin{pf}
Consider~\ref{appprel2e1}. From Lemma~\ref{appprel1}\ref{appprel1e1},
\ref{appprel1e3} follows $\|[\Op]_\um^{-1}\|_s\leq4d^3\Opw_m^{-1}$ and
$\Opw_m^{-1}\leq d\|[\Op]_\um^{-1}\|_s$ which in turn imply
$d^{-1}\leq\|[\Op]_\um^{-1}\|_s\Opw_m \leq D$ and $d^{-1} \leq\Opw_M
\max_{1\leq m\leq M}\|[\Op]_\um^{-1}\|_s\leq D$ due to the monotonicity
of $\Opw$. From these estimates we conclude~\ref{appprel2e1}. Consider
\ref{appprel2e2}. Observe that
$\Delta_{M^+_n}^\Opw\leq\hw_{(M^+_n)}\Opw_{M^+_n }^{-1}$. In case
$M^+_n=1$ the assertion follows from $\hw_{(1)}\Opw_1^{-1}=1$
(Assumption~\ref{assreg}). Thus, let \mbox{$M^{\hw}_n\geq M^+_n>1$}, then
$\min_{1\leq j\leq M^+_n}\{ \Opw_j(j\hw_{(j)})^{-1}\}\geq(1+\log
n)(4Dn)^{-1}$, and hence $M^+_n\Delta_{M^+_n}^\Opw\leq4Dn(1+\log
n)^{-1}$, $\Lambda_{M^+_n }^\Opw \leq(1+\log D)(1+\log n)$,
$M^+_n\Delta_{M^+_n}^{[\Op]}\leq 4D^2n(1+\log n)^{-1}$ and
$\Lambda_{M^+_n}^{[\Op]}\leq(1+2\log D)(1+\log n)$.~\ref{appprel2e2}
follows now by combination of these estimates. Consider
\ref{appprel2e3}. By employing $ D\Opw_{M^+_n }^{-1} \geq\break\max_{1\leq
m\leq M^+_n}\|[\Op]_\um^{-1}\|$, \ref {appprel2e3} follows from
$\Opw_1=1$ if $M^+_n=1$, while for $M^+_n>1$, we use $M^+_n\hw_{(M^+_n
)}\Opw_{M^+_n}^{-1}\leq4Dn(1+\log n)^{-1}$. Consider \ref {appprel2e4}.
Since $\epsilon$ and $X$ are centred the identity
$[\So^{m}]_\um=[\Op]_{\um}^{-1}[g]_{\um}$ implies $\rho_m^2\leq2 (\Ex
Y^2+ \Ex|\langle\So^{m},X\rangle_\Hz|^2 ) = 2 (\sigma_Y^2+[g]_{\um
}^t[\Op ]_{\um}^{-1}[g]_{\um} ) =\sigma_m^2$. By applying successively
the inequality\break $\|\Op^{1/2}\So\|^2\leq d\|\So\|_\Opw^2$ due to \cite
{Heinz1951}, Assumption~\ref{assreg}, that is, $\Opw$ and $b^{-1}$ are
nonincreasing, and the identity $\sigma_Y^2=\sigma^2+\langle
\Op\sol,\sol\rangle_\Hz$ follows
%
%eB.1 #&#
\begin{equation}\label{appprel3e1}
\sigma_Y^2\leq\sigma^2+d\|\sol
\|_\Opw^2\leq \sigma^2+dr.
\end{equation}
Furthermore, from~\ref{appprel1e3} and~\ref{appprel1e5}
in Lemma~\ref{appprel1}, we obtain
%
%eB.2 #&#
\begin{equation}\label{appprel3e2}
[g]_{\um}^t[\Op]_{\um}^{-1}[g]_{\um}
\leq d \bigl\|\So^{m}\bigr\|_\Opw^2
\leq34d^9r,
\end{equation}
which together with (\ref{appprel3e1}) implies~\ref{appprel2e4}
and completes the proof.
\end{pf}
%
%leB.3 #&#
\begin{lem}\label{appprel4}Let $\Op\in\cG_{\Opw}^d$ with $\Opw$
as in
Assumption~\ref{assreg}. For all $n,m\geq1$ holds
\[
{ \biggl\lbrace\frac{1}{4}< \frac{\|[\hOp]_\um^{-1}\|_s}{\|[\Op]_\um^{-1}\|_s}\leq4,\forall 1\leq m\leq
M^{\hw
}_n \biggr\rbrace}\subset{ \bigl\lbrace
M^-_n\leq\Me\leq M^+_n \bigr\rbrace}.
\]
\end{lem}
\begin{pf}
Let $\htau_m:=\|[\hOp]_\um^{-1}\|_s^{-1}$ and
$\tau_m:=\|[\Op]_\um^{-1}\|_s^{-1}$. We use below without
further reference that $D^{-1}\leq\tau_m/\Opw_m \leq d$ due to
Lemma~\ref{appprel2}\ref{appprel2e1}.
The result of the lemma follows by combination of the next two assertions,
%
%eB.3 #&#
%eB.4 #&#
\begin{eqnarray}
\label{appprel4e1}
\bigl\{\Me< M^-_n\bigr\}&\subset& \biggl\{
\min_{1\leq m\leq M^{\hw
}_n}\frac
{\htau_{m}}{\tau_{m}}< \frac{1}{4} \biggr\},
\\
\label{appprel4e2} \bigl\{\Me> M^+_n\bigr\}&\subset& \biggl\{
\max_{1\leq m\leq M^{\hw}_n}\frac
{\htau_{m}}{\tau_{m}}\geq4 \biggr\}.
\end{eqnarray}
Consider (\ref{appprel4e1}) which holds trivially true for $M^-_n=1$.
If $M^-_n>1$, then
\[
\min_{1\leq m\leq
M^-_n} \frac{\Opw_m}{m\hw_{(m)}}\geq\frac{4D(1+\log n)}{n}
\]
implies
\[
\min_{1\leq m\leq M^-_n} \frac{\tau_m}{m\hw_{(m)}}\geq \frac
{4(1+\log n)}{n}
\]
and
\begin{eqnarray*}
{ \bigl\lbrace\Me< M^{\hw}_n \bigr\rbrace}\cap{ \bigl
\lbrace\Me < M^-_n \bigr\rbrace}&=&\bigcup_{M=1}^{M^-_n-1}{
\lbrace\Me =M \rbrace}
\\
&\subset&\bigcup_{M=1}^{M^-_n-1}{ \biggl\lbrace
\frac{\htau_{M+1}}{(M+1)\hw_{(M+1)}}<\frac{1+\log n}{n} \biggr\rbrace}\\
&\subset&{ \biggl\lbrace
\min_{1\leq m\leq M^-_n}\frac{\htau_{m}}{\tau_{m}}< 1/4 \biggr\rbrace},
\end{eqnarray*}
while ${ \lbrace\Me= M^{\hw}_n \rbrace}\cap{
\lbrace\Me< M^-_n \rbrace}=\varnothing$ which shows
(\ref{appprel4e1}) because $M^-_n\leq M^{\hw}_n$.

Consider (\ref{appprel4e2}) which holds trivially true for $M^+_n
=M^{\hw}_n
$. If $M^+_n<M^{\hw}_n$, then
$\frac{\tau_{M^+_n+1}}{(M^+_n+1)\hw_{(M^+_n+1)}}< \frac{(1+\log
n)}{4n}$ and (\ref{appprel4e2}) follows from
\begin{eqnarray*}
{ \lbrace\Me>1 \rbrace}\cap{ \bigl\lbrace\Me> M^+_n \bigr\rbrace}&=&
\bigcup_{M=M^+_n+1}^{M^{\hw}_n}{ \lbrace \Me=M \rbrace}
\\
&\subset&\bigcup_{M=M^+_n+1}^{M^{\hw}_n}{ \biggl\lbrace
\min_{2\leq
m\leq M} \frac{\htau_{m}}{m\hw_{(m)}}\geq\frac{1+\log n}{n}
\biggr\rbrace}\\
&\subset&
{\biggl\lbrace\frac{\htau_{M^+_n+1}}{\tau_{M^+_n+1}}\geq 4 \biggr\rbrace}
\end{eqnarray*}
and $\{\Me= 1\}\cap\{\Me> M^+_n\}=\varnothing$ which
completes the proof.
\end{pf}
%
%leB.4 #&#
\begin{lem}\label{appprel5}
Let $\cA_n$, $\cB_n$ and $\cC_n$ as in (\ref{apppree4}). For all
$n\geq1$ it holds true that $\cA_n\cap\cB_n\cap\cC_n\subset\{
\pen_k\leq\hpen_k\leq72\pen_k,1\leq k\leq M^{\hw}_n\} \cap \{M^-_n
\leq\Me\leq M^+_n\}$.
\end{lem}
\begin{pf}
Let $M^{\hw}_n\geq k\geq1$. If $\|[\Xi]_\uk\|_s\leq1/8$,
that is,
on the event $\cB_n$, it is easily verified that $ \|(\Id_\uk
+[\Xi]_\uk)^{-1}-\Id_\uk\|_s\leq1/7$ which we exploit to conclude
%
%eB.5 #&#
\begin{eqnarray}
\label{appprel5e1}
&\displaystyle \frac{6}{7} \leq\frac{\|[\Diag_{\hw}]_{\uk}^{1/2}[\hOp
]^{-1}_\uk[\Diag_{\hw}]_{\uk}^{1/2}\|_s}{\|[\Diag_{\hw
}]_{\uk}^{1/2}[\Op]^{-1}_\uk[\Diag_{\hw}]_{\uk}^{1/2}\|_s}\leq
\frac{8}{7},\qquad \frac{6}{7} \leq\frac{\|[\hOp]^{-1}_\uk\|_s}{\|[\Op
]^{-1}_\uk\|_s}\leq
\frac{8}{7} \quad\mbox{and}&
\nonumber\\[-8pt]\\[-8pt]
&\displaystyle {6} x^t[\Op]^{-1}_\uk x\leq 7x^t[
\hOp]^{-1}_\uk x \leq8 x^t[\Op
]^{-1}_\uk x\qquad \mbox{for all } x\in\Rz^k&
\nonumber
\end{eqnarray}
and, consequently
%
%eB.6 #&#
\begin{equation}
\label{appprel5e2} (6/7) [\hgf]^t_\uk[
\Op]^{-1}_\uk[\hgf]_\uk\leq[
\hgf]^t_\uk [\hOp ]^{-1}_\uk[
\hgf]_\uk\leq(8/7) [\hgf]^t_\uk[
\Op]^{-1}_\uk[\hgf ]_\uk.
\end{equation}
Moreover, from $\|[\Xi]_\uk\|_{s}\leq1/8$ we obtain after
some algebra,
\begin{eqnarray*}
[g]^t_\uk[\Op]^{-1}_\uk[g]_\uk
&\leq&({1}/{16}) [g]^t_\uk [\Op ]^{-1}_\uk[g]_\uk+
4[W]_\uk[\Op]_\uk^{-1}[W]_{\uk} + 2 [
\hgf ]^t_\uk[\Op]^{-1}_\uk[
\hgf]_\uk,
\\
{}[\hgf]^t_\uk[\Op]^{-1}_\uk[
\hgf]_\uk&\leq&({33}/{16}) [g ]^t_\uk [
\Op]^{-1}_\uk[g]_\uk+ 4[W]_\uk[
\Op]_\uk^{-1}[W]_{\uk}.
\end{eqnarray*}
Combining each of these estimates with (\ref{appprel5e2}) yields
\begin{eqnarray*}
(15/16) [g]^t_\uk[\Op]^{-1}_\uk[g]_\uk
&\leq& 4[W]_\uk[\Op]_\uk^{-1}[W]_{\uk} +
(7/3)[\hgf]^t_\uk[\hOp]^{-1}_\uk[
\hgf]_\uk,
\\
(7/8)[\hgf]^t_\uk[\hOp]^{-1}_\uk[
\hgf]_\uk&\leq&(33/16) [g ]^t_\uk[\Op
]^{-1}_\uk[g]_\uk+ 4[W]_\uk[
\Op]_\uk^{-1}[W]_{\uk}.
\end{eqnarray*}
If in addition $[W]_\uk^t[\Op]_{\uk}^{-1}[W]_\uk\leq\frac
{1}{8}([g]_\uk^t[\Op]^{-1}_\uk[g]_\uk+ \sigma_Y^2) $,
that is, on the event $\csetn$, then the last two
estimates imply, respectively,
\begin{eqnarray*}
(7/16) \bigl([g]^t_\uk[\Op]^{-1}_\uk[g]_\uk+
\sigma_Y^2\bigr) &\leq& (15/16)\sigma_Y^2
+ (7/3) [\hgf]^t_\uk[\hOp]^{-1}_\uk[
\hgf]_\uk,
\\
(7/8)[\hgf]^t_\uk[\hOp]^{-1}_\uk[
\hgf]_\uk&\leq&(41/16) [g ]^t_\uk[\Op
]^{-1}_\uk[g]_\uk+ (1/2)\sigma_Y^2
\end{eqnarray*}
and hence in case ${1}/{2}\leq\hsigma_Y^2/\sigma_Y^2\leq{3}/{2}$,
that is, on the event $\asetn$, we obtain
\begin{eqnarray*}
(7/16) \bigl( [g]^t_\uk[\Op]^{-1}_\uk[g]_\uk+
\sigma_Y^2\bigr)&\leq&(15/8) \hsigma_Y^2
+ (7/3) [\hgf]^t_\uk[\hOp]^{-1}_\uk[
\hgf]_\uk,
\\
(7/8) \bigl([\hgf]^t_\uk[\hOp]^{-1}_\uk[
\hgf]_\uk+ \hsigma{}^2_Y\bigr)&\leq&(41/16)
[g]^t_\uk[\Op]^{-1}_\uk[g]_\uk+
(29/16)\sigma_Y^2.
\end{eqnarray*}
Combining the last two estimates we have
\begin{eqnarray*}
\tfrac{1}{6} \bigl(2[g]^t_\uk[\Op]^{-1}_\uk[g]_\uk+2
\sigma_Y^2\bigr)&\leq& \bigl(2[\hgf]^t_\uk[
\hOp]^{-1}_\uk[\hgf]_\uk+ 2\hsigma{}^2_Y
\bigr)\\
&\leq&3 \bigl(2[g ]^t_\uk[\Op]^{-1}_\uk[g]_\uk+2
\sigma_Y^2\bigr).
\end{eqnarray*}
On $\cA_n\cap\cB_n\cap\cC_n$ the last estimate and (\ref
{appprel5e1}) hold for all $1\leq k\leq M^{\hw}_n$, hence
\[
\cA_n\cap\cB_n\cap\cC_n\subset \biggl\{
\frac{1}{6}\leq\frac
{\hsigma_m^2}{\sigma_m^2}\leq3 \mbox{ and } \frac{6}{7} \leq
\frac{\Delta^{[\widehat{\Op}]}_m}{\Delta_m^{[\Op]}}\leq\frac{8}{7}, \forall1\leq m\leq M^{\hw}_n
\biggr\}.
\]
Moreover it is easily seen that $(6/7) \leq\Delta_m^{[\widehat{\Op
}]}/ \Delta_m^{[\Op]}\leq(8/7)$ implies
\[
1/2\leq\bigl(1+\log(7/6)\bigr)^{-1}\leq\Lambda_m^{[\widehat{\Op}]}/
\Lambda_m^{[\Op]} \leq \bigl(1+\log(8/7)\bigr)\leq3/2.
\]
Due to the last estimates the definitions of $\pen_m$ and $\hpen_m$ imply
%
%eB.7 #&#
\begin{equation}
\label{appprel5e3} \cA_n\cap\cB_n\cap\cC_n
\subset{ \bigl\lbrace\pen_m \leq\hpen_m\leq72
\pen_m, \forall1\leq m\leq M^{\hw}_n \bigr
\rbrace}.
\end{equation}
On the other hand, by exploiting successively (\ref{appprel5e1}) and
Lemma~\ref{appprel4}, we have
\[
\cA_n\cap\cB_n\cap\cC_n\subset{ \biggl
\lbrace\frac{6}{7}\leq \frac{\|[\hOp]^{-1}_\um\|_s}{\|[\Op]^{-1}_\um\|_s}\leq\frac{8}{7}, \forall1
\leq m\leq M^{\hw}_n \biggr\rbrace }\subset{ \bigl\lbrace
M^-_n\leq\Me\leq M^+_n \bigr\rbrace}.
\]
The last display and (\ref{appprel5e3}) imply the assertion of the lemma.
\end{pf}
%
%leB.5 #&#
\begin{lem}\label{appprel6}
For all $m,n\geq1$ with $ n\geq (8/7) \|[\Op]^{-1}_{\um
}\|_s$
we have $\mho_{m,n} \subset\Omega_{m,n}$.
\end{lem}
\begin{pf}
Taking into account
$[\hOp]_{\um}= [\Op]^{1/2}_{\um}\{\Id_{\um}+[\Xi]_{\um}\}[\Op
]^{1/2}_{\um}$ observe that
$\|[\Xi]_{\um}\|_s\leq 1/8$ and $ n\geq(8/7) \|
[\Op]^{-1}_{\um}\|_s$ imply $\|[\hOp]^{-1}_{\um}\|_s
\leq n$ due\vspace*{1pt} to a
Neumann series argument. Hence,
$\mho_{m,n} \subset\Omega_{m,n}$ which proves the lemma.%
\end{pf}

%sC #&#
\section{\texorpdfstring{Proof of Proposition \lowercase{\protect\ref{gaussp1}}}{Proof of Proposition 3.1}}\label{appgauss}
We will suppose throughout this section that the conditions of
Proposition~\ref{gaussp1} are satisfied
which allow us to employ Lemmas~\ref{appprel1}--\ref{appprel6}.
First, we show technical assertions (Lemmas~\ref{appnormal}--\ref
{appgaussl4}) exploiting Assumption~\ref{assgauss}, that is, $X$
and $\epsilon$ are jointly normally distributed. They are used below to
prove that Assumptions~\ref{assremainder-i} and \ref
{assremainder-ii} are satisfied (Propositions~\ref{appgaussp1}
and~\ref{appgaussp2}, resp.), which is the claim of Proposition~\ref
{gaussp1}.

We begin by recalling elementary properties due to Assumption \ref
{assgauss} which are frequently used in this section. Given $f\in
\Hspace$ the random
variable $\langle f,X\rangle_\Hz$ is normally distributed with mean
zero and
variance $\langle\Op f,f\rangle_\Hz$. Consider the Galerkin solution
$\So^{m}
$ and $h\in\Hspace_m$; then
$\langle\sol-\So^{m},X\rangle_\Hz$ and $\langle h,X\rangle_\Hz$
are independent.
Thereby, $Z_m=Y-\langle\So^{m},X\rangle_\Hz=\sigma\epsilon+
\langle\sol-\So^{m},X\rangle_\Hz$ and $[X]_\um$ are independent,
normally distributed with mean zero and, respectively, variance $
\rho_m^2$ and covariance matrix~$[\Op]_\um$. Consequently,
$(\rho_m^{-1}Z_m,[X]_m^t[\Op]_m^{-1/2})$ is a
vector with independent, standard normally distributed entries. The
next assertion states elementary inequalities for Gaussian random
variables and its straightforward proof is omitted.
%
%leC.1 #&#
\begin{lem}\label{appnormal}Let $\{U_i,V_{ij},1\leq i\leq n,1\leq
j\leq m\}$
be independent and standard normally distributed. For all $\eta>0$
and $\zeta\geq4m/n$ we have:

\renewcommand\thelonglist{(\roman{longlist})}
\renewcommand\labellonglist{\thelonglist}
\begin{longlist}
\item\label{appnormaleq1}$P{ ( n^{-1/2}\sum_{i=1}^n
(U_i^2-1)\geq\eta ) }\leq\exp (-\frac{1}{8}\frac{\eta^2}{1+\eta
n^{-1/2}} )$;
\item\label{appnormaleq2}$P{ ( n^{-1}\llvert \sum_{i=1}^n U_iV_{i1}\rrvert \geq\eta ) }\leq\frac{\eta
n^{1/2}+1}{\eta n^{1/2} } \exp
(-\frac{n}{4}\min{ \lbrace\eta^2,1/4 \rbrace} )$;
\item\label{appnormaleq2b}$P (n^{-2}\sum_{j=1}^m\llvert
\sum_{i=1}^n U_iV_{ij}\rrvert^2\geq\zeta )\leq\exp
(\frac
{-n}{16} )
+\exp (\frac{-\zeta n}{64} )$;
\end{longlist}
and for all $c\geq1$ and $a_1,\ldots,a_m\geq0$ we obtain:

\renewcommand\thelonglist{(\roman{longlist})}
\renewcommand\labellonglist{\thelonglist}
\begin{longlist}
\setcounter{longlist}{3}
\item\label{appnormaleq3}$\Ex (\sum_{i=1}^n U_i^2-2 c
n
)_+\leq16\exp (\frac{-c n}{16} )$;\vadjust{\goodbreak}
\item\label{appnormaleq4}$\Ex (\sum_{j=1}^m\llvert
n^{-1/2}\sum_{i=1}^n U_iV_{ij}\rrvert^2- 4  c  m)_+\leq
16\exp (\frac{-c
m}{16} )
+ 32 \frac{c  m}{n} \exp (\frac{-n}{16} )$;
\item\label{appnormaleq5}$\Ex ( \sum_{j=1}^m a_j\llvert
\sum_{i=1}^n U_iV_{ij}\rrvert^2 )^2=n(n+2) (\sum_{j=1}^ma_j^2 +
(\sum_{j=1}^ma_j)^2 )$.\vspace*{-1pt}
\end{longlist}
\end{lem}
%
%leC.2 #&#
\begin{lem}\label{appgaussl1} For all $n,m\geq1$ we have:

\renewcommand\thelonglist{(\roman{longlist})}
\renewcommand\labellonglist{\thelonglist}
\begin{longlist}
\item\label{appgaussl1e1}$n^2\rho_m^{-4}\Ex\|[W]_\um\|^4\leq
6  (\Ex\| X\|^2)^2$.
\end{longlist}
Furthermore, there exist a numerical constant $C>0$ such that for all
$n\geq1$:

\renewcommand\thelonglist{(\roman{longlist})}
\renewcommand\labellonglist{\thelonglist}
\begin{longlist}
\setcounter{longlist}{1}
\item\label{appgaussl1e2}
$n^8\max_{1\leq m\leq{\lfloor n^{1/4}\rfloor}} P (\frac{[W]_\um^t[\Op]_\um^{-1}[W]_\um}{\rho_m^2}>\frac{1}{16} )\leq C$;
\item\label{appgaussl1e3}
$ n^8 \max_{1\leq m\leq{\lfloor n^{1/4}\rfloor}} P{ ( \|
[\Xi]_\um\|_s> 1/8 ) }\leq C$;
\item\label{appgaussl1e4}
$n^7P{ ( \{1/2\leq\hsigma{}^2_Y/\sigma^2_Y\leq3/2\}^c )
}\leq C$.\vspace*{-1pt}
\end{longlist}
\end{lem}
\begin{pf}
Denote by $(\lambda_j,e_j)_{1\leq j\leq m}$ an eigenvalue decomposition
of $[\Op]_\um$. Define $U_i:=(\sigma\epsilon_i +
\langle\sol-\So^{m},X_i\rangle_\Hz)/\rho_m$ and
$V_{ij}:=(\lambda_j^{-1/2}e_j^t[X_i]_{\um})$, $1\leq i\leq n$, $1\leq
j\leq m$, where $U_1,\ldots,U_n,V_{11},\ldots,V_{nm}$ are independent
and standard normally distributed. Consider~\ref{appgaussl1e1} and \ref
{appgaussl1e2}. Taking into account $\sum_{j=1}^m\lambda_j\leq\Ex \|
X\|_\Hz^2$ and the identities $n^4\rho_m^{-4}\|[W]_\um\|^4 =
(\sum_{j=1}^m \lambda_j(\sum_{i=1}^nU_iV_{ij})^2)^2$ and
$([W]_\um^t[\Op]_\um^{-1}[W]_\um)/\rho_m^2 =n^{-2} \sum_{j=1}^m
(\sum_{i=1}^nU_iV_{ij})^2$, assertions~\ref{appgaussl1e1} and
\ref{appgaussl1e2} follow, respectively, from in Lemma
\ref{appnormal}\ref{appnormaleq5} and~\ref{appnormaleq2b} (with
$a_j=\lambda_j$). Consider~\ref{appgaussl1e3}. Since
$n\|[\Xi]_\um\|_{s}\leq m\max_{1\leq j,l\leq
m}|\sum_{i=1}^n(V_{ij}V_{il}-\delta_{jl})|$ we obtain due to
\ref{appnormaleq1} and~\ref{appnormaleq2} in Lemma \ref
{appnormal} that for all $\eta>0$%
\begin{eqnarray*}
&&
P\bigl(\bigl\|[\Xi]_\um\bigr\|_{s}\geq\eta\bigr)\\[-1pt]
&&\quad\leq\sum
_{1\leq j,l\leq
m}P\Biggl(\Biggl|n^{-1}\sum
_{i=1}^n(V_{ij}V_{il}-
\delta_{jl})\Biggr|\geq\eta/m\Biggr)
\\[-1pt]
&&\quad\leq m^2\max \Biggl\{ P\Biggl(\Biggl|\frac{1}{n}\sum
_{i=1}^nV_{i1}V_{i2}\Biggr|\geq
\frac{\eta}{m}\Biggr), P\Biggl(\Biggl|\frac
{1}{n^{1/2}}\sum
_{i=1}^n\bigl(V_{i1}^2-1
\bigr)\Biggr|\geq n^{1/2}\frac{\eta
}{m}\Biggr) \Biggr\}
\\[-1pt]
&&\quad\leq m^2\max \biggl\{\biggl(1+ \frac{m}{\eta n^{1/2}}\biggr) \exp
\biggl(-\frac{n}{4}\min\biggl\{\frac{\eta^2}{m^2},\frac{1}{4}\biggr
\} \biggr), 2 \exp \biggl(-\frac{1}{8}\frac{n\eta^2/m^2}{1+\eta/m} \biggr) \biggr
\}.%
\end{eqnarray*}
Keeping in mind that $1/8=\eta\leq m/2$, the last bound implies \ref
{appgaussl1e3}.
Consider~\ref{appgaussl1e4}. Since $\{Y_i/\sigma_Y\}_{i=1}^n$ are
independent, standard, normally
distributed and
$\{1/2\leq\hsigma{}^2_Y/\sigma^2_Y\leq3/2\}^c\subset\{|n^{-1}\sum_{i=1}^nY_i^2/\sigma_Y^2-1|>1/2\}$,~\ref{appgaussl1e4} follows
from Lem\-ma~\ref{appnormal}\ref{appnormaleq1}.\vspace*{-1pt}
\end{pf}
%
%leC.3 #&#
\begin{lem}\label{appgaussl2} We have for all $c\geq1$ and $n,m\geq1$
\[
\Ex \biggl( \frac{n[W]_\um^t[\Op]_\um^{-1}[W]_\um
}{\rho_m^2} - 4 c m \biggr)_+\leq16\exp \biggl(
\frac{-c m}{16} \biggr) + 32 \frac
{c
m}{n} \exp \biggl(\frac{-n}{16}
\biggr).\vspace*{-1pt}
\]
\end{lem}

\begin{pf}
From\vspace*{1pt}
$n\|[\Op]_\um^{-1/2}[W]_\um\|^2\rho_m^{-2} = \sum_{j=1}^m
( n^{-1/2}\sum_{i=1}^n U_iV_{ij})^2$ derived in the proof of Lemma
\ref{appgaussl1} and Lemma~\ref{appnormal}\ref{appnormaleq4} follows the
assertion.\vadjust{\goodbreak}
\end{pf}
%
%leC.4 #&#
\begin{lem}\label{appgaussl3} There is a constant $C(d)$ depending
on $d$ such that for all $n\geq1$,
\[
\sup_{\sol\in\mathcal{F}^r_b}\sup_{\Op\in\cG_{\Opw}^d}\sum_{
k={m^\diamond_n}}^{M^+_n
}
\Delta_{k}^{[\Op]} \Ex \biggl([W]_{\uk}^t[
\Op]_{\uk}^{-1}[W]_{\uk} - 4 \sigma_k^2
\frac{k \Lambda_k^{[\Op]}}{n} \biggr)_+ \leq C(d) \bigl(\sigma^2+r\bigr) \Sigma
n^{-1}.
\]
\end{lem}
\begin{pf}
The key argument of the proof is Lemma~\ref{appgaussl2} with
$c=\Lambda_k^{[\Op]}$. Taking into account this bound and for all
$\So\in\Soclass_{\Sow }^r$ and $\Op \in \Opclass_{\Opw}^d$ that
$\Delta_k^{[\Op]}\leq4d^3 \Delta_k^\Opw$, $ (1+\log d)^{-1}
\Lambda_k^\Opw\leq\Lambda_k^{[\Op]}$, $ \delta_{M^+_n}^{[\Op ]}\leq n C
d^6(1+\log d)$ and $\rho_k^2\leq\sigma_{k}^2\leq2(\sigma^2+35d^6r)$
[Lemma~\ref{appprel2}(i), (ii) and (iv), resp.] hold true, we obtain
\begin{eqnarray*}
&&
\sum_{ k={m^\diamond_n}}^{M^+_n}\Delta_{k}^{[\Op]}
\Ex \biggl([W]_{\uk
}^t[\Op ]_{\uk}^{-1}[W]_{\uk}
- 4 \sigma_k^2 \frac{k \Lambda_k^{[\Op]
}}{n} \biggr)_+
\\[-1pt]
&&\qquad\leq\sum_{ k=1}^{M^+_n}\frac{\sigma_k^2\Delta_{k}^{[\Op]}}{n} \Ex
\biggl(\frac{n[W]_{\uk}^t[\Op]_{\uk}^{-1}[W]_{\uk}}{\rho_k^2} - 4 k \Lambda_k^{[\Op]} \biggr)_+
\\[-1pt]
&&\qquad\leq C(d) \bigl(\sigma^2+r\bigr)\\[-1pt]
&&\qquad\quad{}\times n^{-1} \Biggl\{ \sum
_{ k=1}^{M^+_n} \Delta_k^\Opw
\exp \biggl(-\frac
{k\Lambda_k^\Opw}{16(1+\log d)} \biggr) + M^+_n\exp ({-n}/{16} )
\Biggr\}.
\end{eqnarray*}
Finally, exploiting the constant $\Sigma$ satisfying (\ref{defSigma})
and $M^+_n\exp ({-n}/{16} )\leq C$ for all $n\geq1$, we obtain
the assertion of the lemma.
\end{pf}
%
%leC.5 #&#
\begin{lem}\label{appgaussl4} There exist a numerical constant $C$
and a constant $C(d)$ only depending on $d$ such that for all
$n\geq1$, we have:

\renewcommand\thelonglist{(\roman{longlist})}
\renewcommand\labellonglist{\thelonglist}
\begin{longlist}
\item\label{appgaussl4e1}$\sup_{\sol\in\mathcal{F}^r_b}\sup_{\Op\in
\cG_{\Opw}^d} { \lbrace n^6 (M^+_n)^2 \max_{1\leq m\leq M^+_n}
P{ ( \mho_{m,n}^c ) } \rbrace}\leq C$;
\item\label{appgaussl4e2}$\sup_{\sol\in\mathcal{F}^r_b}\sup_{\Op\in
\cG_{\Opw}^d}{ \lbrace n M^+_n\max_{1\leq m\leq M^+_n} P{
( \Omega_{m,n}^c ) } \rbrace}\leq
C(d)$;
\item\label{appgaussl4e3}$\sup_{\sol\in\mathcal{F}^r_b}\sup_{\Op\in
\cG_{\Opw}^d}{ \lbrace n^7 P{ ( \esetn^c ) }
\rbrace}\leq C$.
\end{longlist}
\end{lem}
\begin{pf}
Since $M^+_n\leq{\lfloor n^{1/4}\rfloor}$ and $\mho_{m,n}^c={
\lbrace\|[\Xi]_{\um}\|>1/8 \rbrace}$ assertion \ref {appgaussl4e1}
follows from Lemma~\ref{appgaussl1}\ref{appgaussl1e2}. Consider
\ref{appgaussl4e2}. Let $ n_o:=n_o(d):=\exp(128 d^6)\geq8d^3$, and
consequently $\hw_{(M^+_n)}(M^+_n\log n)\geq128d^6$ for all $n\geq
n_o$. We distinguish in the following the cases $n<n_o$ and $n\geq
n_o$. First, let $1\leq n\leq n_o$. Obviously, $M^+_n\max_{1\leq m\leq
M^+_n} P(\Omega^c_{m,n})\leq M^+_n\leq n^{-1} n_o^{5/4}\leq C(d)n^{-1}$
since $M^+_n\leq n^{1/4}$ and $n_o$ depends on $d$ only. On the other
hand, if $n\geq n_o$, then Lemma~\ref{appprel2}(iii) implies
$n\geq2\max_{1\leq m\leq M^+_n}\|[\Op]^{-1}_\um\|$, and hence
$\mho_{m,n}\subset\Omega_{m,n}$ for all $1\leq m\leq M^+_n$ by
employing Lemma~\ref{appprel6}. From~\ref{appgaussl4e1} we
conclude\break
$M^+_n\max_{1\leq m\leq M^+_n}P(\Omega^c_{m,n})\leq M^+_n\max_{1\leq
m\leq M^+_n }P(\mho^c_{m,n})\leq C n^{-3}$. By\vspace*{1pt} combination of the two
cases we obtain~\ref{appgaussl4e2}. It remains to show
\ref{appgaussl4e3}. Consider the events $\cA_n$, $\cB_n$ and $\cC_n$
given in~(\ref{apppree4}), where $\cA_n\cap\cB_n\cap\cC_n\subset\cE_n$
due to Lemma~\ref{appprel5}. We have $n^7 P{ ( \cA^c_n ) }\leq C$, $n^7
P{ ( \cB^c_n ) }\leq C$, $n^7 P{ ( \cC^c_n ) }\leq C$ due to Lemma
\ref{appgaussl1}\ref{appgaussl1e4},~\ref{appgaussl1e3},
\ref{appgaussl1e2}, respectively (keep in mind ${\lfloor
n^{1/4}\rfloor}\geq M^{\hw}_n$ and
$2(\sigma_Y^2+[g]_\uk^t[\Op]_\uk^{-1}[g]_\uk)=\sigma_{k}^2\geq\rho_k^2$).
Combining these estimates implies~\ref{appgaussl4e3}.
\end{pf}
%
%prC.6 #&#
\begin{prop}\label{appgaussp1}Let $\kappa=96$ in definition (\ref
{defpen}) of the penalty $\pen$. There exists a constant $C(d)$ only
depending on $d$ such that for all $n\geq1$,
\[
\sup_{\sol\in\mathcal{F}^r_b}\sup_{\Op\in\cG_{\Opw}^d} \Ex { \biggl\lbrace\sup_{{m^\diamond_n} \leq k\leq M^+_n}
\biggl(\bigl\| \hsol_{k}-\So^{k}\bigr\|^2_\hw-
\frac{1}{6}\pen_k \biggr)_+ \biggr\rbrace}\leq C(d) \bigl(
\sigma^2+r\bigr) \Sigma n^{-1}.
\]
\end{prop}
\begin{pf}
Since $[\hsol_{k}-\So^{k}]_\uk=
[\hOp]^{-1}_\uk[W]_{\uk}\one_{\Omega_{k,n}} - [\So^{k} ]_\uk\one_{\Omega_{k,n}^c}$ it follows
%
%eC.1 #&#
\begin{equation}
\label{appgaussp1decomp} \bigl\|\hsol_k-\So^{k}
\bigr\|_\hw^2= \bigl\|[\Diag_\hw]_\uk^{1/2}[
\hOp]^{-1}_\uk[W]_{\uk}\bigr\|^2
\one_{\Omega_{k,n}} +\bigl\| \So^{k} \bigr\|^2_\hw
\one_{\Omega_{k,n}^c}.
\end{equation}
Exploiting $\|(\Id_{\uk}+[\Xi]_{\uk})^{-1}\|_s\one_{\mho _{k,n}}\leq
2$, $[\hOp]_{\uk}= [\Op]^{1/2}_{\uk}\{\Id_{\uk}+[\Xi]_{\uk}\} [\Op
]^{1/2}_{\uk}$ and the definition of $\Delta_k^{[\Op]}$ imply $
\|[\Diag_\hw
]^{1/2}_\uk[\hOp]^{-1}_\uk[W]_{\uk}\|^2\one_{\mho_{k,n}}\leq4
\Delta_k^{[\Op]} \|[\Op]^{-1/2}_\uk[W]_{\uk}\|^2$. On the other hand,
we have $
\|[\Diag_\hw]^{1/2}_\uk[\hOp]^{-1}_\uk[W]_{\uk}\|^2\one_{\Omega_{k,n}}\leq\hw_{(k)}
n^2 \|[W]_{\uk}\|^2$. From these estimates and $\|\So^{k}
\|_\hw\leq\|\So^{k}\|_b$ ($\hw b^{-1}$ is nonincreasing due to
Assumption~\ref{assreg}) we deduce for all $k\geq1$,
\[
\bigl\|\hsol_k-\So^{k}\bigr\|_\hw^2
\leq4\Delta_k^{[\Op]}\bigl\|[\Op ]_{\uk}^{-1/2}[W]_{\uk}
\bigr\|^2+\hw_{(k)} n^2\bigl\|[W]_{\uk}
\bigr\|^2\one_{\mho
_{k,n}^c} + \bigl\|\So^{k}
\bigr\|^2_b\one_{\Omega_{k,n}^c}.
\]
This upper bound
and $\pen_k=96\sigma_k^2 k\Delta_k^{[\Op]}\Lambda_k^{[\Op]}
n^{-1}$ imply
\begin{eqnarray*}
&&
\Ex \biggl\{\sup_{{m^\diamond_n}\leq k\leq M^+_n} \biggl(\bigl\|\hsol_{k}-
\So^{k} \bigr\|^2_\hw-\frac{\pen_k}{6} \biggr)_+
\biggr\} \\
&&\qquad\leq\sum_{k={m^\diamond_n}}^{M^+_n}
n^3 \bigl(\Ex\bigl\|[W]_{\uk
}\bigr\|^4
\bigr)^{1/2} \bigl(P\bigl(\mho_{k,n}^c\bigr)
\bigr)^{1/2}
+\sum_{k={m^\diamond_n}}^{M^+_n} \bigl\|\So^{k}
\bigr\|^2_bP\bigl(\Omega^c_{k,n}
\bigr)\\
&&\qquad\quad{} +4\sum_{ k={m^\diamond_n}}^{M^+_n}
\Delta_{k}^{[\Op]} \Ex \biggl(\bigl\|[\Op]_{\uk}^{-1/2}[W]_{\uk}
\bigr\|^2 - 4 \sigma_k^2 \frac
{k \Lambda_k^{[\Op]}}{n}
\biggr)_+.
\end{eqnarray*}
By exploiting Lemmas~\ref{appprel1}\ref{appprel1e5} and
\ref{appgaussl1}\ref{appgaussl1e1} together with $ \rho_m^2\leq
2(\sigma^2+35d^6r)$ [Lemma~\ref{appprel2}(iv)] the first and
second r.h.s. term are bounded by
\begin{eqnarray*}
&&
6\bigl(\sigma^2+35d^6r\bigr)\Ex\| X
\|^2 n^2M^+_n\max_{{m^\diamond
_n}\leq
k\leq M^+_n
} \bigl(P
\bigl(\mho_{k,n}^c\bigr) \bigr)^{1/2}\\
&&\qquad{} +
34d^8rM^+_n\max_{{m^\diamond_n}
\leq
k\leq M^+_n} P\bigl(
\Omega^c_{k,n}\bigr).\vadjust{\goodbreak}
\end{eqnarray*}
Combining this upper bound, the property $\Ex\| X\|^2\leq d
\sum_{j\geq1}\Opw_j\leq d\Sigma$ and the estimates given in Lemma
\ref
{appgaussl4}, we deduce for all $\sol\in\mathcal{F}^r_b$ and $\Op
\in\cG_{\Opw}^d
$ that
\begin{eqnarray*}
&&
\sup_{\sol\in\mathcal{F}^r_b}\sup_{\Op\in\cG_{\Opw}^d} \Ex { \biggl\lbrace\sup_{{m^\diamond_n} \leq k\leq M^+_n}
{ \biggl( \bigl\| \hsol_{k}-\So^{k}\bigr\|^2_\hw-
\frac{1}{6}\pen_k \biggr)_{+}} \biggr\rbrace} \\
&&\qquad\leq
C(d) \bigl(\sigma^2+r\bigr) \Sigma n^{-1}
\\
&&\qquad\quad{}+ 4 \sup_{\sol\in\mathcal{F}^r_b}\sup_{\Op\in\cG_{\Opw}^d}\sum_{ k={m^\diamond_n}
}^{M^+_n}
\Delta_{k}^{[\Op]} \Ex{ \biggl( \bigl\|[\Op]_{\uk}^{-1/2}[W]_{\uk}
\bigr\|^2 - 4 \sigma_k^2 \frac{k \Lambda_k^{[\Op]}}{n}
\biggr)_{+}}.
\end{eqnarray*}
The result of the proposition follows now by replacing the last r.h.s.
term by its upper bound given in Lemma~\ref{appgaussl3}, which
completes the proof.
\end{pf}
%
%prC.7 #&#
\begin{prop}\label{appgaussp2} Let $\kappa=96$ in definition (\ref
{defpen}) and (\ref{defhpenMen}) of $\pen$ and $\hpen$. There exists
a constant $C(d)$ only depending on $d$ such that for all
$n\geq1$
\[
\sup_{\sol\in\mathcal{F}^r_b}\sup_{\Op\in\cG_{\Opw}^d} \Ex \bigl(\|\hsol_{\whm}-
\sol\|_\hw^2\one_{\cE^c_n}\bigr) \leq C(d) \bigl(
\sigma^2+r\bigr) \Sigma n^{-1}.
\]
\end{prop}
\begin{pf}
From the decomposition (\ref{appgaussp1decomp}) and $ \|[\Diag_\hw]^{1/2}_\uk[\hOp]^{-1}_\uk[W]_{\uk}\|^2\one_{\Omega
_{k,n}}\leq\Delta^\hw_k
n^2 \|[W]_{\uk}\|^2$ given in the proof of Proposition \ref
{appgaussp1} we conclude
\[
\|\hsol_k-\sol\|_\hw^2\leq2
\Delta^\hw_k n^2 \bigl\| [W]_{\uk}
\bigr\|^2 + 2 \bigl\|\So^{k}\bigr\|^2_\hw+
2 \|\sol\|^2_\hw\qquad \mbox{for all } k\geq1.
\]
By exploiting Lemma~\ref{appprel1}\ref{appprel1e5} together with
$\|\So^{k}\|_\hw\leq\|\So^{k}\|_b$ ($\hw b^{-1}$ is
nonincreasing due to Assumption
\ref{assreg}) we obtain for all $\sol\in\mathcal{F}^r_b$ and $\Op
\in\cG_{\Opw}^d
$%
\[
\|\hsol_k-\sol\|_\hw^2\leq2
\Delta^\hw_k n^2 \bigl\| [W]_{\uk}
\bigr\|^2 + 2 \bigl(34 d^8r+r\bigr)\qquad \mbox{for all } k
\geq1.
\]
Since $1\leq\whm\leq M^{\hw}_n$ and $\max_{1\leq k\leq M^{\hw}_n}
\hw_{(k)}\leq
n$ it follows that
\begin{eqnarray*}
\Ex\bigl(\|\hsol_{\whm}-\sol\|_\hw^2
\one_{\cE^c_n}\bigr)&\leq&2 n^3 M^{\hw}_n
\max_{1\leq k\leq M^{\hw}_n} \bigl(\Ex\bigl\|[W]_{\uk}\bigr\|^4
\bigr)^{1/2}\bigl|P\bigl(\cE^c_n\bigr)\bigr|^{1/2}\\
&&{}
+ 70 d^8rM^{\hw}_nP\bigl(\cE^c_n
\bigr).
\end{eqnarray*}
From Lemma~\ref{appgaussl1}\ref{appgaussl1e1} together with $
\rho_m^2\leq2(\sigma^2+35d^6r)$ (Lemma~\ref{appprel2}) and
$\Ex\| X\|^2\leq d\Sigma$, we conclude for all $\sol\in
\mathcal{F}^r_b$
and $\Op\in\cG_{\Opw}^d$ that
\begin{eqnarray*}
\Ex\bigl(\|\hsol_{\whm}-\sol\|_\hw^2
\one_{\cE^c_n}\bigr)&\leq&12 \bigl(\sigma^2+35d^6r
\bigr) \,d\Sigma n^2 M^{\hw}_n\bigl|P\bigl(
\cE^c_n\bigr)\bigr|^{1/2} \\
&&{}+ 70 d^8r
M^{\hw}_nP\bigl(\cE^c_n\bigr).
\end{eqnarray*}
The result of the proposition follows now from $M^{\hw}_n\leq{\lfloor
n^{1/4}\rfloor} $ and by replacing the probability $P(\cE^c_n)$ by
its upper
bound $Cn^{-7}$ given in Lemma~\ref{appgaussl4}.
\end{pf}
\begin{pf*}{Proof of Proposition~\ref{gaussp1}}
The assertion follows from Propositions~\ref{appgaussp1} and
\ref{appgaussp2},
and we omit the details.\vadjust{\goodbreak}
\end{pf*}
%
%sD #&#
\section{\texorpdfstring{Proof of Proposition \lowercase{\protect\ref{momp1}}}{Proof of Proposition 3.3}}\label{appmom}
We assume throughout this section that the conditions of Proposition
\ref{momp1} are satisfied which allows us
to employ Lemmas~\ref{appprel1}--\ref{appprel6}. We formulate first
preliminary results (Proposition~\ref{appgeneraltala} and Lemmas
\ref{appgenerall1}--\ref{appgenerall4}) relying on the moment
conditions (Assumption~\ref{assmom}). They are used to prove that
Assumptions
\ref{assremainder-i} and~\ref{assremainder-ii} are satisfied
(Propositions~\ref{appgeneralp1}
and~\ref{appgeneralp2}, resp.), which is the claim of Proposition
\ref
{momp1}. We begin by gathering elementary bounds due to Assumption
\ref
{assmom}. Let $k$ be given by Assumption~\ref{assmom}; then for all
$m\geq1$ we have
\begin{eqnarray*}
&\displaystyle \Ex|Z_m|^{4k}\leq\rho_m^2
\eta^{4k},\qquad \Ex|Y|^{4k}\leq\sigma^{4k}_Y
\eta^{4k},&\\
&\displaystyle \max_{1\leq j\leq m}\Ex \bigl|\bigl([\Op]_\um^{-1/2}[X]_\um
\bigr)_j\bigr|^{4k}\leq\eta^{4k},&
\\
&\displaystyle \Ex \bigl|\bigl\langle\sol-\So^{m},X\bigr\rangle_\Hz
\bigr|^{4k}\leq\bigl\| \Op^{1/2}\bigl(\So^{m}-\sol\bigr)
\bigr\|_\Hz^{4k} \eta^{4k},&\\
&\displaystyle  \Ex
\bigl|[X]_\um^t[\Op]_\um^{-1}[X]_\um
\bigr|^{2k}\leq m^{2k}\eta^{4k}.&
\end{eqnarray*}
From $\Ex|V|\one_{\{|V|\geq t\}}\leq t^{-k+1}\Ex|V|^k$, $t>0$, under
Assumption~\ref{assmom} follows%
%
%eD.1 #&#
\begin{eqnarray}
\label{appgenerale1}
&\displaystyle \Ex\epsilon^2\one_{\{|\epsilon|>n^{1/6}\}}\leq
\frac{\eta^{32}}{n^{5}},&
\nonumber
\\
&\displaystyle \Ex \bigl|\bigl\langle\sol-\So^{m},X\bigr\rangle_\Hz
\bigr|^{2}\one_{\{
|\langle\sol-\So^{m},X\rangle_\Hz|>\|\Op^{1/2}(\So^{m}-\sol
)\|_\Hz n^{1/6}\}} &\nonumber\\[-8pt]\\[-8pt]
&\qquad\qquad\displaystyle \leq \frac
{\eta^{32}}{n^5}\bigl\|
\Op^{1/2}\bigl(\So^{m}-\sol\bigr)\bigr\|_\Hz^2,\hspace*{110pt}&
\nonumber\\
&\displaystyle \Ex\bigl|[X]_\um^t[\Op]_\um^{-1}[X]_\um\bigr|^2
\one_{\{[X]_\um
^t[\Op]_\um^{-1}[X]_\um>m n^{1/3}\}}
\leq\frac{\eta^{32}}{n^{14/3}}m^{2}&
\nonumber
\end{eqnarray}
for all $m,n\geq1$, and by employing Markov's inequality
%
%eD.2 #&#
\begin{eqnarray}
\label{appgenerale2} P\bigl(|\epsilon|>n^{1/6}\bigr)&\leq&
\frac{\eta^{32}}{n^{16/3}},\nonumber\\[-8pt]\\[-8pt]
P\bigl(\bigl|\bigl\langle \sol-\So^{m},X\bigr
\rangle_\Hz\bigr|>\bigl\|\Op^{1/2}\bigl(\So^{m}-\sol
\bigr)\bigr\|_\Hz n^{1/6}\bigr)&\leq& \frac
{\eta^{32}}{n^{16/3}}.\nonumber
\end{eqnarray}
We exploit these bounds in the following proofs. The key argument used
in the proof of Lemma~\ref{appgenerall2} is the following inequality
due to~\cite{Talagrand1996}; see, for example,~\cite{KleinRio}.
%
%prD.1 #&#
\begin{prop}[(Talagrand's inequality)]
\label{appgeneraltala}
Let $T_1,\ldots, T_n$ be independent $\cT$-valued random variables
and $\nu^*_s
= (1/n)\sum_{i=1}^n [\nu_s(T_i) - \Ex[\nu_s(T_i)]  ]$, for
$\nu_s$
belonging to a countable class $\{\nu_s\dvtx s\in\Sz\}$ of measurable
functions. Then,
for $\epsilon> 0$,
\begin{eqnarray*}
&&
\Ex{ \Bigl( \sup_{s\in\Sz} \bigl|\nu^*_s\bigr|^2 - 2(1+2
\epsilon)H^2 \Bigr) }_+
\\
&&\qquad\leq C \biggl(\frac{v}{n}\exp\biggl( -K_1\epsilon
\frac{nH^2}v \biggr) + \frac{h^2}{n^2C^2(\epsilon)} \exp\biggl(-K_2C(
\epsilon)\sqrt{\epsilon}\frac{nH}{h}\biggr) \biggr)
\end{eqnarray*}
with $K_1 = 1/6$, $K_2= 1/(21\sqrt{2})$, $C(\epsilon) =
\sqrt{1+\epsilon}-1$ and $C$ a universal constant and where
\[
\sup_{s\in\Sz}\sup_{t\in\cT}\bigl|\nu_s(t)\bigr| \leq h,\qquad \Ex
\Bigl[\sup_{s\in\Sz}\bigl|\nu^*_s\bigr| \Bigr]\leq H,\qquad
\sup_{s\in
\Sz} \frac{1}{n}\sum_{i=1}^n
\Var\bigl(\nu_s(T_i)\bigr)\leq v.
\]
\end{prop}
%
%leD.2 #&#
\begin{lem}\label{appgenerall1}There exist a numerical constant $C>0$
such that for all $n\geq1$:

\renewcommand\thelonglist{(\roman{longlist})}
\renewcommand\labellonglist{\thelonglist}
\begin{longlist}
\item\label{appgenerall1e1}$n^2\sup_{m\geq1}\rho_m^{-4}\Ex
\|[W]_\um\|^4\leq C \eta^8 (\Ex\| X\|_\Hz^2)^2$;
\item\label{appgenerall1e2}$n^8\max_{1\leq m\leq{\lfloor
n^{1/4}\rfloor}}
P{ ( \frac{[W]_\um^t[\Op]_\um^{-1}[W]_\um}{\rho_m^2}>\frac
{1}{16} ) }\leq C \eta^{64}$;
\item\label{appgenerall1e3}$n^8 \max_{1\leq m\leq{\lfloor
n^{1/4}\rfloor}}
P{ ( \|[\Xi]_\um\|_s> 1/8 ) }\leq C(\eta)$;
\item\label{appgenerall1e4}$n^7P{ ( \{1/2\leq\hsigma{}^2_Y/\sigma^2_Y\leq3/2\}^c ) }\leq C \eta^{64}$.
\end{longlist}
\end{lem}
\begin{pf}
Denote by $(\lambda_j,e_j)_{1\leq j\leq m}$ an eigenvalue decomposition
of $[\Op]_\um$. Define $U_i:=(\sigma\epsilon_i +
\langle\sol-\So^{m},X_i\rangle_\Hz)/\rho_m$ and $V_{ij}:=(\lambda_j^{-1/2}e_j^t[X_i]_{\um})$, $1\leq i\leq n$, $1\leq j\leq m$. Keep in
mind that $\Ex|U_i|^{4k}\leq\eta^{4k}$, $\Ex|V_{ij}|^{4k}\leq\eta^{4k}$ and $\Ex|U_iV_{ij}|^{2k}\leq\eta^{4k}$
for $k\geq16$ (Assumption~\ref{assmom}) and $\{U_iV_{ij}\}_{i=1}^n$
are independent, centred for $1\leq j\leq m$. Consider
\ref{appgenerall1e1},~\ref{appgenerall1e2} where $n^4\rho_m^{-4}\|[W]_\um\|^4= (\sum_{j=1}^m\lambda_j(\sum_{i=1}^nU_iV_{ij})^2)^2$ and
$([W]_\um^t[\Op]_\um^{-1}[W]_\um)/\rho_m^2 =n^{-2} \sum_{j=1}^m
(\sum_{i=1}^nU_iV_{ij})^2$. Applying\vspace*{2pt} Minkowski's (resp., Jensen's)
inequality and Theorem 2.10 in~\cite{Petrov1995}, we have
\begin{eqnarray*}
{n^2}\rho^{-4}_m\Ex\bigl\|[W]_\um
\bigr\|^4 &\leq& n^{-2} \Biggl[\sum
_{j=1}^m \lambda_j \Biggl(\Ex \Biggl|\sum
_{i=1}^nU_iV_{ij}
\Biggr|^4 \Biggr)^{1/2} \Biggr]^2 \\
&\leq& C
\eta^8 \Biggl[\sum_{j=1}^m
\lambda_j \Biggr]^2;
\\
{n^{k}} {m^{-k}\rho^{-2k}_m}\Ex\bigl\|[
\Op]_\um^{-1/2}[W]_\um\bigr\|^{2k} &\leq&
\frac{1}{ m}\sum_{j=1}^m
n^{-k}\Ex \Biggl|\sum_{i=1}^nU_iV_{ij}
\Biggr|^{2k} \leq C(k) \eta^{4k},
\end{eqnarray*}
which, respectively, implies~\ref{appgenerall1e1}, since $\sum_{j=1}^m\lambda_j\leq\Ex\| X\|_\Hz^2$, and~\ref{appgenerall1e2},
by employing Markov's inequality.
Proof of~\ref{appgenerall1e3}. Since $\{V_{ij}V_{il}-\delta_{jl}\}_{i=1}^n$ are
independent, centred with $\Ex|V_{ij}V_{il}-\delta_{jl}|^{2k}\leq C
\eta^{4k}$, $1\leq j,l\leq m$, Theorem 2.10 in~\cite{Petrov1995}
implies $n^{k}\Ex |n^{-1}\sum_{i=1}^n(V_{ij}V_{il}-\delta_{jl})
|^{2k}\leq C(k)\eta^{4k}$ and $m^{-2k} n^{k}\Ex\|[\Xi]_\um
\|_s^{2k}\leq C(k)
\eta^{4k}$ because $\|[\Xi]_m\|_s^2\leq\sum_{1\leq j,l\leq
m}|V_{ij}V_{il}-\delta_{jl}|^2$.
Applying Markov's inequality gives~\ref{appgenerall1e3}.
Proof of~\ref{appgenerall1e4}. Since $\{Y_i^2/\sigma_Y^{2}-1\}_{i=1}^n$ are independent, centred with
$\Ex |Y_i^2/\sigma_Y^2-1 |^{2k}\leq C(k) \eta^{4k}$ Theorem
2.10 in
\cite{Petrov1995} implies\break $\Ex |n^{-1}\sum_{i=1}^nY_i^2/\sigma_Y^2-1 |^{2k}\leq C(k) n^{-k} \eta^{4k}$ and
$P (|{n^{-1}}\sum_{i=1}^nY_i^2/\sigma_Y^2-1|>\break 1/2 )\leq C n^{-16}
\eta^{64}$ employing Markov's inequality.~\ref{appgaussl1e4} follows
now from $\{1/2\leq\hsigma{}^2_Y/\sigma^2_Y\leq3/2\}^c\subset\{
|{n^{-1}}\sum_{i=1}^nY_i^2/\sigma_Y^2-1|>1/2\}$.
\end{pf}
%
%leD.3 #&#
\begin{lem}\label{appgenerall2} Let $\varsigma_m:= \sigma+\eta^2\|\Op^{1/2}(\So^{m}-\sol)\|_\Hz$, $m\geq1$. There exists a
numerical constant $C$ such that
for all ${\lfloor n^{1/4}\rfloor}\geq m\geq1$ we have
\begin{eqnarray*}
&&
\Ex \biggl(\frac{\|[\Op]_{\um}^{-1/2}[W_n]_{\um
}\|^2}{\varsigma_m^2} -12\frac{m\Lambda_m^{[\Op]}}{n} \biggr)_+\\
&&\qquad\leq
\frac
{C}{n} \biggl\{ \exp \biggl( - \frac{m\Lambda_m^{[\Op]}}{6} \biggr) + \exp
\biggl( - \frac{n^{1/6}}{100} \biggr) + \frac{\eta^{32} }{n^{2}} \biggr\}.
\end{eqnarray*}
\end{lem}
\begin{pf}
Let $\Sz^m:=\{z\in\Rz^m\dvtx z^tz\leq1\}$. Define $\cE_n:=\{e\in\Rz
\dvtx|e|\leq
n^{1/6}\}$, $\cX_{1n}:=\{x\in\Hz\dvtx |\langle\sol-\So^{m},x\rangle_\Hz|\leq
\|\Op^{1/2}(\sol-\So^{m})\|_\Hz n^{1/6}\}$, $\cX_{2n}:=\{
x\in
\Hz\dvtx\break [x]_\um^t[\Op]_\um^{-1}[x]_\um\leq m n^{1/3}\}$ and
$\cX_n:=\cX_{1n}\cap\cX_{2n}$. For $e\in\Rz$, $x\in\Hz$, $s\in
\Sz^m$
set%
\begin{eqnarray*}
\nu_{s}(e,x)&:=& \bigl(\sigma e+ \bigl\langle\sol-
\So^{m},x\bigr\rangle_\Hz\bigr) s^t[
\Op]_\um^{-1/2}[x]_\um\one_{\{e\in\cE_n,x\in\cX_n\}},
\\
R_{s}(e,x)&:=& \bigl(\sigma e+ \bigl\langle\sol-\So^{m},x
\bigr\rangle_\Hz\bigr) s^t[\Op]_\um^{-1/2}[x]_\um(1-
\one_{\{e\in\cE_n,x\in\cX_n\}}).
\end{eqnarray*}
Let $\nu^*_s:= n^{-1}\sum_{i=1}^n\{\nu_{s}(\epsilon_i,X_i) - \Ex
\nu_{s}(\epsilon_i,X_i)\}$ and $R^*_s:= n^{-1}\sum_{i=1}^n\{
R_{s}(\epsilon_i,X_i) - \Ex R_{s}(\epsilon_i,X_i)\}$, then $\|[\Op]_{\um
}^{-1/2}[W_n]_{\um}\|^2=\sup_{s\in\Sz^m}|\nu_{s}^*
+R_{s}^*|^2$ and hence
%
%eD.3 #&#
\begin{eqnarray}
\label{appgenerall2decomp}
&&
\Ex \biggl(\bigl\|[\Op]_{\um}^{-1/2}[W_n]_{\um}
\bigr\|^2 -12\varsigma_m^2\frac{
m\Lambda_m^{[\Op]}}{n}
\biggr)_+ \nonumber\\
&&\qquad\leq2\Ex \biggl( \sup_{s\in\Sz^m}\bigl|\nu_{s}^*
\bigr|^2 -6\varsigma_m^2\frac{
m\Lambda_m^{[\Op]}}{n}
\biggr)_+
+ 2\Ex\sup_{s\in\Sz
^m}\bigl|R_{s}^*\bigr|^2\\
&&\qquad=:2\{
T_1+T_2\},
\nonumber
\end{eqnarray}
where we bound the r.h.s. terms $T_1$ and $T_2$ separately. Consider
first $T_1$. We intend to apply Talagrand's inequality. To this end,
for $e\in\Rz,x\in\Hz$, we have%
%
%eD.4 #&#
\begin{eqnarray}\label{appgenerall2e11}\qquad
\sup_{s\in\Sz^m} \bigl|\nu_s(e,x)\bigr|^2 &=& \bigl(\sigma e+
\bigl\langle\sol-\So^{m},x\bigr\rangle_\Hz
\bigr)^2[x]_\um^t[\Op ]_\um^{-1}[x]_\um
\one_{\{e\in\cE_n,x\in\cX_{n}\}}
\nonumber\\[-8pt]\\[-8pt]
&\leq&\bigl(\sigma+ \bigl\|\Op^{1/2}\bigl(\So^{m}-\sol\bigr)
\bigr\|_\Hz\bigr)^2 n^{2/3} m \leq
\varsigma^2_m n^{2/3} m=:h^2.
\nonumber
\end{eqnarray}
By employing the independence of $\epsilon$ and $X$ it is easily seen that
\begin{eqnarray*}
n \Ex\sup_{ s\in\Sz^m} \bigl|\nu_{s}^*\bigr|^2&\leq&
\sigma^2m +\Ex \bigl|\bigl\langle\sol-\So^{m},X\bigr
\rangle_\Hz\bigr|^2 [X]^t_{\um}[
\Op]_\um^{-1}[X]_\um,
\\
\sup_{ s\in\Sz^m} \frac{1}{n} \sum_{i=1}^n
\Var\bigl(\nu_{s}(\epsilon_i,X_i)\bigr)
&\leq&\sigma^2 + \sup_{ s\in\Sz^m} \Ex\bigl|\bigl\langle\sol-
\So^{m},X\bigr\rangle_\Hz\bigr|^2\bigl|s^t[
\Op]_\um^{-1/2}[X]_\um\bigr|^2.
\end{eqnarray*}
By applying the Cauchy--Schwarz inequality together with $\Ex\|
[\Op]_\um^{-1/2}[X]_\um\|^{4}\leq m^2\eta^{4}$ and $\Ex
|\langle\sol-\So^{m},X\rangle_\Hz |^{4}\leq\|\Op^{1/2}(\So^{m}-\sol)\|_\Hz^{4} \eta^{4}$
we obtain
%
%eD.5 #&#
\begin{equation}
\label{appgenerall2e12}\quad \Ex\sup_{ s\in\Sz^m} \bigl|\nu_{s}^*\bigr|^2
\leq\frac{m}{n} \bigl(\sigma^2 + \bigl\|\Op^{1/2}\bigl(
\sol-\So^{m}\bigr)\bigr\|_\Hz^2
\eta^4\bigr)\leq\varsigma_m^2
\frac{m
\Lambda_m^{[\Op]}}{n}=:H^2,
\end{equation}
and taking into account that $\Ex|s^t[\Op]_\um^{-1/2}[X]_\um
|^{4}\leq\eta^{4}$, $ s\in\Sz^m$, we obtain
%
%eD.6 #&#
\begin{equation}
\label{appgenerall2e13}\quad \sup_{ s\in\Sz^m} \frac{1}{n} \sum
_{i=1}^n \Var\bigl(\nu_{s}(
\epsilon_i,X_i)\bigr) \leq\sigma^2 + \bigl\|
\Op^{1/2}\bigl(\So^{m}-\sol\bigr)\bigr\|_\Hz^{2}
\eta^4\leq\varsigma_m^2=:v.
\end{equation}
Due to (\ref{appgenerall2e11})--(\ref{appgenerall2e13})
Talagrand's inequality (Lemma~\ref{appgeneraltala} with $\epsilon
=1$) implies
%
%eD.7 #&#
\begin{equation}
\label{appgenerall2e2} \Ex \biggl(\sup_{ s\in\Sz^m} \bigl|\nu_{s}^*\bigr|^2-6
\varsigma^2_m \frac
{m\Lambda_m^{[\Op]}}{n} \biggr)_+ \leq C
\frac{\varsigma_m^2}{n} \biggl\{\exp \biggl( - \frac{m\Lambda_m^{[\Op]
}}{6} \biggr) + \exp
\biggl( - \frac{n^{1/6}}{100} \biggr) \biggr\},\hspace*{-35pt}
\end{equation}
where we used that $m\leq{\lfloor n^{1/4}\rfloor}$. Consider $T_2$.
By employing $[X]_\um[\Op]_\um^{-1}[X]_\um\*\one_{\{X\in\cX
_{2,n}\}
}\leq
mn^{1/3}$ and $\cX_{n}=\cX_{1n}\cap\cX_{2n}$ we have
\begin{eqnarray*}
n\Ex\sup_{ s\in\Sz^m} \bigl|R_{s}^*\bigr|^2 &\leq&\Ex\bigl(\sigma
\epsilon+ \bigl\langle\sol-\So^{m},X\bigr\rangle_\Hz
\bigr)^2 [X]_\um [\Op]_\um^{-1}[X]_\um
\one_{\{X\notin\cX_{2,n}\}}
\\
&&{}+ m n^{1/3}\Ex\bigl(\sigma\epsilon+ \bigl\langle\sol-
\So^{m},X\bigr\rangle_\Hz\bigr)^2(
\one_{\{\epsilon\notin\cE
_n\}
}+\one_{\{
X\notin\cX_{1n}\}}).
\end{eqnarray*}
Since $ \Ex(\sigma\epsilon+ \langle\sol-\So^{m},X\rangle_\Hz
)^4\leq
(\sigma^2+\|\Op^{1/2}(\sol-\So^{m})\|_\Hz^2)^2\eta^4$, $\Ex
\epsilon^2=1$ and
$\Ex
|\langle\sol-\So^{m},X\rangle_\Hz|^2=\|\Op^{1/2}(\sol-\So^{m} )\|_\Hz^2$ the
independence of $\epsilon$ and $X$ implies
\begin{eqnarray*}
n\Ex\sup_{ s\in\Sz^m} \bigl|R_{s}^*\bigr|^2&\leq&\bigl(
\sigma^2+\bigl\|\Op^{1/2}\bigl(\sol-\So^{m}\bigr)
\bigr\|_\Hz^2\bigr)\eta^2 \bigl(
\Ex\bigl|[X]_\um[\Op ]_\um^{-1}[X]_\um\bigr|^2
\one_{\{X\notin\cX_{2,n}\}} \bigr)^{1/2}
\\
&&{}+ m n^{1/3} \bigl\{ \sigma^2 \Ex\epsilon^2
\one_{\{\epsilon\notin
\cE_n\}
} +\bigl\|\Op^{1/2}\bigl(\sol-\So^{m}\bigr)
\bigr\|_\Hz^2P(\epsilon\notin\cE_n)
\\
&&\hspace*{42.8pt}{}+ \sigma^2 P(X\notin\cX_{1n})+ \Ex\bigl|\bigl\langle\sol-
\So^{m},X\bigr\rangle_\Hz\bigr|^2
\one_{\{X\notin\cX_{1n}\}} \bigr\}.
\end{eqnarray*}
We exploit now the estimates given in (\ref{appgenerale1}) and
(\ref{appgenerale2}). Thereby, we obtain
\[
n\Ex\sup_{ s\in\Sz^m} \bigl|R_{s}^*\bigr|^2 \leq C \bigl(
\sigma^2+\bigl\|\Op^{1/2}\bigl(\sol-\So^{m}\bigr)
\bigr\|_\Hz^2\bigr)\eta^{32} m n^{-7/3}
\leq C \varsigma^2_m\eta^{32} n^{-2},
\]
where we used that $m\leq{\lfloor n^{1/4}\rfloor}$. Keeping in mind
decomposition (\ref{appgenerall2decomp}), the last bound and (\ref
{appgenerall2e2})
imply together the claim of Lemma~\ref{appgenerall2}.
\end{pf}
%
%leD.4 #&#
\begin{lem}\label{appgenerall3} There exists a constant $K:=K(\sigma
,\eta,\mathcal{F}^r_b,\cG_{\Opw}^d)$ depending on $\sigma$, $\eta
$ and the classes
$\mathcal{F}^r_b$ and $\cG_{\Opw}^d$ only
such that for all $n\geq1$ we have
\[
\sup_{\sol\in\mathcal{F}^r_b}\sup_{\Op\in\cG
_{\Opw}^d}\sum
_{
m={m^\diamond_n}
}^{M^+_n} \Delta_m^{[\Op]}\Ex
\biggl(\bigl\|[\Op]_{\um
}^{-1/2}[W_n]_{\um}
\bigr\|^2 -\frac{12\sigma_m^2 m\Lambda_m^{[\Op]}}{n} \biggr)_+\leq K\frac{\eta^{32}
(\sigma^2+r)
\Sigma}{n}.
\]
\end{lem}
\begin{pf}
There exists an integer $n_o:=n_o(\sigma,\eta,\mathcal{F}^r_b,\cG_{\Opw}^d)$ depending
on $\sigma$, $\eta$ and the classes $\mathcal{F}^r_b$ and $\cG_{\Opw}^d$ only such that
for all $n\geq n_o$ and for all $m\geq{m^\diamond_n}$ we have
$\varsigma_m^2\leq2(\sigma^2 + \|\Op^{1/2}\sol\|_\Hz^2 +[g
]_{\um
}^t[\Op]_{\um}^{-1}[g]_{\um})= 2(\sigma^2_Y +[g]_{\um}^t[\Op
]_{\um
}^{-1}[g]_{\um})=\sigma_m^2$. Indeed, we
have $1/{m^\diamond_n}=o(1)$ as $n\to\infty$ and $|\varsigma_m^2-\sigma^2|=o(1)$ as $m\to\infty$ because $\varsigma_m= \sigma+\eta^2\|\Op^{1/2}(\So^{m}-\sol)\|_\Hz$ and $\|\Op^{1/2}(\So^{m}-\sol)\|_\Hz^2\leq
34  d^9   r
\Opw_mb_m^{-1}$ due to Lemma~\ref{appprel1}\ref{appprel1e5}.
% n_o$.
First, consider $n<n_o$. Due to Lemma
\ref{appgenerall1}\ref{appgenerall1e1} and $
\rho_m^2\leq2(\sigma^2+35d^6r)$ [Lemma~\ref{appprel2}(iv)] we have for
all $m\geq1$
\[
\Ex \biggl(\bigl\|[\Op]_{\um}^{-1/2}[W_n]_{\um}
\bigr\|^2 -\frac{
12\sigma_m^2m\Lambda_m^\Op}{n} \biggr)_+\leq\Ex\bigl\|[
\Op]_{\um
}^{-1/2}[W_n]_{\um}
\bigr\|^2\leq C \frac{m}{n} \eta^{4}\bigl(
\sigma^2+d^6r\bigr).
\]
Hence, $M^+_n\leq{\lfloor n^{1/4}\rfloor}$ and
$m\Delta_m^\Op\leq\delta_{M^+_n}^\Op\leq n C(d)$ [Lemma \ref
{appprel2}(ii)] imply
\begin{eqnarray*}
&&\sup_{\sol\in\mathcal{F}^r_b}\sup_{\Op\in\cG_{\Opw}^d}\sum_{
m={m^\diamond_n}}^{M^+_n
}
\Delta_m^{[\Op]}\Ex \biggl(\bigl\|[\Op]_{\um}^{-1/2}[W_n]_{\um}
\bigr\|^2 -12\sigma_m^2\frac{ m\Lambda_m^{[\Op]}}{n}
\biggr)_+\\
&&\qquad\leq C(d)\frac
{n_o^{5/4}\eta^4(\sigma^2+r)}{n},
\end{eqnarray*}
which proves the lemma for all $1\leq n< n_o$.
Consider now $n\geq n_o$ where $\varsigma_m^2\leq\sigma_m^2$ for all
$m\geq{m^\diamond_n}$. Thereby, we can apply Lemma~\ref{appgenerall2},
which gives
\begin{eqnarray*}
&&
\sup_{\sol\in\mathcal{F}^r_b}\sup_{\Op\in\cG_{\Opw}^d}\sum_{
m={m^\diamond_n}}^{M^+_n
}
\Delta_m^{[\Op]}\Ex \biggl(\bigl\|[\Op]_{\um}^{-1/2}[W_n]_{\um}
\bigr\|^2 -12\sigma_m^2\frac{ m\Lambda_m^{[\Op]}}{n}
\biggr)_+
\\
&&\qquad\leq C\sup_{\sol
\mathcal{F}^r_b
}\sup_{\Op\in\cG_{\Opw}^d}\sum_{ m={m^\diamond_n}}^{M^+_n}
\frac
{\varsigma_m^2\Delta_m^{[\Op]
}}{n} \biggl\{ \exp \biggl( - \frac{m\Lambda_m^{[\Op]}}{6} \biggr) + \exp
\biggl( - \frac{n^{1/6}}{100} \biggr) + \frac{\eta^{32} }{n^{2}} \biggr\}.
\end{eqnarray*}
Since $\Delta_k^{[\Op]}\leq4d^3 \Delta_k^\Opw$, $\Lambda_k^{[\Op]
}\geq(1+\log d)^{-1} \Lambda_k^\Opw$, $M^+_n\Delta_{M^+_n }^{[\Op]}\leq
\delta_{M^+_n}^{[\Op]}\leq n C d^6(1+\log d)$ and
$\varsigma_k^2\leq\sigma_{k}^2\leq2(\sigma^2+35d^6r)$
[Lemma~\ref{appprel2}(i), (ii), (iv), resp.] follows
\begin{eqnarray*}
&&
\sup_{\sol\in\mathcal{F}^r_b}\sup_{\Op\in\cG_{\Opw}^d}\sum_{
m={m^\diamond_n}}^{M^+_n
}
\Delta_m^{[\Op]}\Ex \biggl(\bigl\|[\Op]_{\um}^{-1/2}[W_n]_{\um}
\bigr\|^2 -12\sigma_m^2\frac{ m\Lambda_m^{[\Op]}}{n}
\biggr)_+ \\
&&\qquad\leq C(d) \bigl(\sigma^2+r\bigr) n^{-1}
\\
&&\qquad\quad{}\times\sup_{\sol\in\mathcal{F}^r_b}\sup_{\Op\in\cG_{\Opw}^d} \Biggl\{ \sum
_{
m={m^\diamond_n}
}^{M^+_n} \Delta_m^\Opw\exp
\biggl(-\frac{m\Lambda_m^\Opw}{6(1+\log d
)} \biggr) + n\exp \biggl( - \frac{n^{1/6}}{100} \biggr) +
\frac{\eta^{32} }{n} \Biggr\}.
\end{eqnarray*}
Finally, $\Sigma=\Sigma(\cG_{\Opw}^d)$ as in (\ref{defSigma}) and
$n\exp
(-{n^{1/6}}/{100} )\leq C$ imply for $n\geq n_o$
\begin{eqnarray*}
&&
\sup_{\sol\in\mathcal{F}^r_b}\sup_{\Op\in\cG_{\Opw}^d}\sum_{
m={m^\diamond_n}}^{M^+_n
}
\Delta_m^{[\Op]}\Ex \biggl(\bigl\|[\Op]_{\um}^{-1/2}[W_n]_{\um}
\bigr\|^2 -12\sigma_m^2\frac{ m\Lambda_m^{[\Op]}}{n}
\biggr)_+ \\
&&\qquad\leq C(d)\frac{\eta^{32}(\sigma^2+r)\Sigma}{n}.
\end{eqnarray*}
Combining the cases $n<n_o$ and $n\geq
n_o$ completes the proof.
\end{pf}
%
%leD.5 #&#
\begin{lem}\label{appgenerall4} There exist a numerical constant $C$
and a constant $C(d)$ only depending on $d$ such that for all
$n\geq1$ we have:

\renewcommand\thelonglist{(\roman{longlist})}
\renewcommand\labellonglist{\thelonglist}
\begin{longlist}
\item\label{appgenerall4e1}$\sup_{\sol\in\mathcal{F}^r_b}\sup_{\Op\in
\cG_{\Opw}^d}
\{n^6 (M^+_n)^2 \max_{1\leq m\leq M^+_n} P{ ( \mho_{m,n}^c ) } \}
\leq C\eta^{64}$;
\item\label{appgenerall4e2}$\sup_{\sol\in\mathcal{F}^r_b}\sup_{\Op\in
\cG_{\Opw}^d}
\{n M^+_n\max_{1\leq m\leq M^+_n} P{ ( \Omega_{m,n}^c ) }
\}\leq
C(d)\eta^{64}$;
\item\label{appgenerall4e3}$\sup_{\sol\in\mathcal{F}^r_b}\sup_{\Op\in
\cG_{\Opw}^d}
\{n^7 P{ ( \cE^c_n ) } \}\leq C\eta^{64}$.
\end{longlist}
\end{lem}
\begin{pf}
By employing Lemma~\ref{appgenerall1} rather than Lemma \ref
{appgaussl1} the proof follows along the lines of the proof of Lemma
\ref{appgaussl4}, and we omit the details.
\end{pf}
%
%prD.6 #&#
\begin{prop}\label{appgeneralp1} Let $\kappa=288$ in the definition
(\ref{defpen}) of the penalty $\pen$. There exists a constant
$K:=K(\sigma,\eta,\mathcal{F}^r_b,\cG_{\Opw}^d)$ depending on
$\sigma$, $\eta$
and the
classes $\mathcal{F}^r_b$ and $\cG_{\Opw}^d$ only
such that for all $n\geq1$, we have
\[
\sup_{\sol\in\mathcal{F}^r_b}\sup_{\Op\in\cG_{\Opw}^d} \Ex { \biggl\lbrace\sup_{{m^\diamond_n} \leq k\leq M^+_n}
\biggl(\bigl\| \hsol_{k}-\So^{k}\bigr\|^2_\hw-
\frac{1}{6}\pen_k \biggr)_+ \biggr\rbrace}\leq K
\eta^{64} \bigl(\sigma^2+r\bigr) \Sigma n^{-1}.
\]
\end{prop}
\begin{pf}
By employing Lemmas~\ref{appgenerall1},~\ref{appgenerall3} and
\ref{appgenerall4} rather than Lemmas~\ref{appgaussl1},
\ref{appgaussl3} and~\ref{appgaussl4} the proof follows along the lines
of the proof of Proposition~\ref{appgaussp1}, and we omit the details.
\end{pf}
%
%prD.7 #&#
\begin{prop}\label{appgeneralp2} Let $\kappa=288$ in definition
(\ref{defpen}) and (\ref{defhpenMen}) of $\pen$ and $\hpen$.
There exists
a constant $C(d)$ only depending on $d$ such that for all $n\geq1$,
\[
\sup_{\sol\in\mathcal{F}^r_b}\sup_{\Op\in\cG_{\Opw}^d} \Ex \bigl(\|\hsol_{\whm}-
\sol\|_\hw^2\one_{\cE^c_n}\bigr) \leq C(d)
\eta^{64} \bigl(\sigma^2+r\bigr) \Sigma n^{-1}.
\]
\end{prop}
\begin{pf}
Taking into account Lemmas~\ref{appgenerall1}\ref{appgenerall1e1}
and~\ref{appgenerall4} rather than Lemmas
\ref{appgaussl1}\ref{appgaussl1e1} and~\ref{appgaussl4}
the proof follows along the lines of the proof of Proposition \ref
{appgaussp2}, and we omit the details.\vadjust{\goodbreak}
\end{pf}
\begin{pf*}{Proof of Proposition~\ref{momp1}}
The result follows from Propositions~\ref{appgeneralp1} and~\ref{appgeneralp2}, and we omit the details.
\end{pf*}
\end{appendix}

\section*{Acknowledgments}

We are grateful to two referees and the Associate Editor for
constructive criticism and clear guidelines.

\begin{supplement}%[id=suppA]
\stitle{Simulation study}
\slink[doi]{10.1214/12-AOS1050SUPP} %[doi,text={...}] - jei reikia
%suskaldyti doi
\sdatatype{.pdf}
\sfilename{aos1050\_supp.pdf}
\sdescription{A simulation study illustrating the finite sample
behavior of the fully data-driven estimation procedure and its good
performance.}
\end{supplement}

% imsref loaded by lrinkeviciute, 2012-11-21 14:02:15
% imsref loaded by lrinkeviciute, 2012-11-21 14:48:28

\printaddresses

\end{document}